\input amstex
\documentstyle{amsppt}
\TagsOnRight %\NoBlackBoxes
\loadbold
\define\<{\langle}
\define\>{\rangle}
\define\e{\epsilon}
\define\p{\partial}
\define\pw{\partial_\omega}
 \loadbold
 \define\W{\pmb{|} }
 \define\zd{\Bbb Z^d}
 
 \define\un{\underline}
\let\\=\par

\document
\tenpoint
\topmatter
\title {A  KAM theorem with applications to partial  equations of higher dimension}
\endtitle
\rightheadtext{ A  KAM theorem with application to PDEs}
\author {Xiaoping Yuan}\endauthor
\affil { School of Mathematical Sciences, Fudan University,
Shanghai 200433, China;\\ Key Lab of Mathematics for
Nonlinear Science (Fudan University), Ministry of Education, China.\\
Email:
 xpyuan\@fudan.edu.cn;\ yuanxiaoping\@hotmail.com }\endaffil
\address {Xiaoping Yuan,
School of mathematical Sciences,  Fudan University, Shanghai
200433, China, \\ Key Laboratory of Mathematics for Nonlinear
Science (Fudan University), Ministry of Education.}
\endaddress
%\date Oct.28, 2004\enddate
\email { yuanxiaoping\@hotmail.com;\ xpyuan\@fudan.edu.cn
}\endemail
\thanks Supported by NNSFC  and NCET-04-0365.
\endthanks
\endtopmatter

{\bf Abstract.} {\it The existence of lower dimensional KAM tori
is shown for a class of nearly integrable Hamiltonian systems
where the second Melnikov's conditions are eliminated. As a
consequence, it is proved that there exist many  invariant tori
and thus quasi-periodic solutions for nonlinear wave equations,
Schr\"odinger equations and other equations of higher spatial
dimension. }

\subhead 1. Introduction and main results \endsubhead

Let us begin with the non-linear wave (NLW) equation
$$u_{tt}-u_{xx}+V(x)u+ h(x,u)=0 \tag 1.1$$
subject to Dirichlet  boundary conditions.
  The existence of solutions, periodic in time, for NLW equations
  has been studied by many authors. See [B, B-B, B-G,B-P, Br, L-S] and the references
  theirin, for example. While  finding quasi-periodic solutions, one
  will inevitably
  encounter so-called small divisor difficulty.  The KAM (Kolmogorov-Arnold-Moser) theory is a very powerful
   tool to overcome the difficulty.
This theory deals with the existence of invariant tori ( and thus
quasi-periodic solutions ) for nearly integrable Hamiltonian
systems. In order to obtain the quasi-periodic solutions of a
partial differential equation, one may show the existence of the
lower (finite) dimensional invariant tori for the infinitely
dimensional Hamiltonian system defined by the partial differential
equation. Now consider a Hamiltonian of the form:
$$H=(\omega,y)+\sum_{j=1}^\iota \Omega_jz_j\bar z_j+ R(x,y,z,\bar z)
,\ \iota\le \infty, \tag 1.2$$ with  tangential frequency vector
$\omega=(\omega_1,...,\omega_n)$ and normal frequency vector
$\Omega=(\Omega_1,...,\Omega_\iota)$. When $R\equiv 0$, there is a
trivial invariant torus $x=\omega t, y=0,z=\bar z=0$. The KAM
theory guarantees the persistence of the trivial invariant torus
with a small deformation under a sufficiently small perturbation
$R$, provided that the well-known Melnikov conditions are
fulfilled:
$$( k,\omega)-\Omega_j\neq 0\qquad \text{(the first
Melnikov's)}\tag 1.3$$ for all $k\in\Bbb Z^n$ and $1\le
j\le\iota$, and
$$( k,\omega)+ \Omega_{j_1}- \Omega_{j_2}\neq 0\qquad
\text{(the second Melnikov's)}\tag 1.4$$
 for all $k\in\Bbb Z^n$ and $1\le j_1,j_2\le\iota, j_1\neq j_2$.
 See [E,K1,P1,W] for the details.
This KAM theorem can be applied to a wide array of Hamiltonian
partial
 differential equations of $1$-dimensional spatial variable,
 including (1.1). In that line, Kuksin[K1,2] shows that there are many quasi-periodic solutions
 of (1.1), assuming  that
    the   potential $V$ depends on an $n$-dimensional external
    parameter vector in
    some  non-degenerate way.  Wayne[W] obtains also the existence
    of the quasi-periodic solutions
   of  (1.1), when the
   potential $V$  is lying    on the
   outside of the set of some ``bad" potentials. In [W], the set of all potentials
   is given some Gaussian measure and then the set of ``bad" potentials is
   of   small measure.
    Bobenko \& Kuksin[Bo-K] and P\"oschel[P2]
     get the existence of invariant tori and quasi-periodic
     solutions for a given potential $V(x)\equiv m\in (0,\infty)$. By the remark in [P2], the same result holds also true
   for the parameter values $-1<m<0$.  When
   $m\in(-\infty,-1)\setminus\Bbb Z$, it is shown in [Y1] that
   there are many hyperbolic-elliptic invariant tori. More
   recently, the existence of invariant tori ( thus quasi-periodic
   solutions) of (1.1) are shown for any prescribed non-vanishing and smooth potential
   \footnote{ This potential $V$ contains no parameter.} $V(x)$ in
   [Y2] and for $V(x)\equiv 0$ in [Y3]. In [C-W], [Bo1], [C-Y] and [Br-K-S], the equation (1.1) subject to
   periodic boundary conditions is investigated.

  For NLW equation (1.1) of  spatial dimension $1$, the multiplicity of normal frequency $\Omega_j$ is $1$ in Dirichlet
  boundary condition and $2$ in periodic boundary condition. Considering partial differential
  equations with spatial dimension$>1$, a significant new
  problem arises due to the presence of clusters of normal
  frequencies of the Hamiltonian system defined by partial differential
  equations. In this
  case,  the multiplicity of $\Omega_j$ goes to $\infty$ as $|j|\to\infty$; consequently,
  the second Melnikov's conditions is destroyed seriously, thus preventing the application of the KAM
  theorems mentioned above to  Hamiltonian partial differential equations of higher
  spatial dimension.
  Bourgain[Bo1-4] develops another profound approach, originally proposed by Craig-Wayne in [C-W], and
  successfully obtains the existence of quasi-periodic solutions
  of the nonlinear Schr\"odinger (NLS) equations and NLW equations of
  higher dimension  in space. Instead of KAM theory, the used
  approach is  based on a generalization of  Lyapunov-Schmidt
  procedure and a technique by Fr\"ohlich and
  Spencer[F-S]. Bambusi calls this approach
   ``C-W-B method".

  The advantage of the KAM approach is, on one hand, to possibly
  simplify the proof and, on the other hand, to allow the construction
  of local normal forms closed to the considered torus, which
  could be useful for the better understanding of the dynamics. For
  example, in generally, one can easy check the linear stability
  and the vanishing Lyapunov exponents. As a counterpart of the KAM theory,
  an advantage of the C-W-B
  method is more flexible than the KAM scheme to deal with
  resonant cases where the second Melnikov's conditions are
  violated seriously, although a (local) normal
  form can not be obtained by  using the C-W-B method.
 % And in generally, one does not know whether   or not the quasi-periodic solutions obtained by C-W-B method are on a rotational invariant torus.
  Therefore, it is expected to find out a method by
  which the resonant cases can be easily dealt and by which the
  normal form can be obtained at the same time. In that direction,
  the first result is due to Bourgain[Bo5], to my knowledge. In order to introduce Bourgain's
  idea in [Bo5], let's recall the basic idea of lower dimensional  KAM tori. See [E],
   [K1] and  [P1] for details.  Consider a finitely
  dimensional Hamiltonian (1.2) with $\iota<\infty$. Decompose the perturbation $R$ in
  (1.2) into as follows
  $$\align R=&R^x(x)+(R^y(x),y)\tag 1.5\\&+\<R^z(x),z\>+\<R^{\bar z}(x),\bar
  z\>\tag 1.6\\&+\<R^{zz}(x)z,z\>+
  \<R^{z\bar z}(x)z,\bar z\>+\<R^{zz}(x)\bar z,\bar z\>\tag 1.7\\&+
  O(|y|^2+|y||z|+|z|^3),\tag 1.8\endalign$$
where $(x,y,z,\bar z)$ is in the some subset of $\Bbb
T^n\times\Bbb R^n\times\Bbb C^\iota\times\Bbb C^\iota$, and where
$R^x, R^y: \Bbb T^n\to\Bbb C^n$ and $R^z, R^{\bar z}:\Bbb
T^n\to\Bbb C^\iota$ and $R^{zz}, R^{z\bar z}, R^{\bar z\bar z}$
are complex $\iota\times\iota$  matrices. In the classic KAM
theory, one finds a series of symplectic transformations to kill
the perturbing terms (1.5), (1.6) and (1.7) such that the
transformed Hamiltonian is of the form
$$\un H=(\omega,y)+\sum_{j=1}^\iota \Omega_jz_j\bar z_j+O(|y|^2+|y||z|+|z|^3).
\tag 1.9 $$ One sees easily that $\Bbb
T^n\times\{y\}|_{y=0}\times\{z,\bar z\}_{z=\bar z=0}$ is an
invariant torus of (1.9). Going back to (1.2), one can get an
invariant torus of (1.2). The second Melnikov's conditions (1.4)
comes out while killing the perturbing term (1.7).  Bourgain[Bo5]
modifies essentially the idea mentioned above. He eliminates the
perturbing terms (1.5) and (1.6) and puts the perturbing term
(1.7) into the ``integrable" part of the Hamiltonian $H$.
Consequently, he obtains a Hamiltonian of the form
$$\align\tilde H=&(\omega,y)+\sum_{j=1}^\iota \Omega_jz_j\bar z_j+\<R^{zz}(x)z,z\>+
  \<R^{z\bar z}(x)z,\bar z\>+\<R^{zz}(x)\bar z,\bar z\>\tag 1.10\\&
+O(|y|^2+|y||z|+|z|^3).\tag 1.11\endalign $$ This Hamiltonian
$\tilde H$ can be regarded as a counterpart of the local normal
form of the classic KAM theory. It is clear that
  $\Bbb T^n\times\{y=0\}\times\{z=\bar z=0\}$ is
  still
  an invariant torus of $\tilde H$. Because of without killing
  (1.7), the second Melnikov's conditions (1.4) are not required.
Therefore, using the Bourgain's idea, one can deal resonant cases
where the second Melnikov's conditions (1.4) are violated
seriously; at the same time, one can get a normal form (1.10+11).
However, the ``integrable" term (1.10) is not really integrable,
since it contains the angle-variable $x$ in $R^{zz}(x)$'s. Because
of this fact, while eliminating (1.6), one will have to solve a
homological equation with variable coefficients
$$\sqrt{-1}(\omega,\p_x)F^z+(\Omega^0
+R^{zz}(x)+\cdots)F^z=R^z \tag 1.12$$ where $F^z$ is unknown
function. By contrast, in the classical KAM theorem, the
homological equation is of constant coefficients. To solve (1.12)
one needs to investigate the inverse of a ``big" matrix $A$ of the
form:
$$A=\text{ diag }((k,\omega)+\Omega_j^0: k\in\Bbb Z^n,j=1,...,\iota)+
(\widehat{R^{zz}}(k-l):k,l\in\Bbb Z^n).\tag 1.13$$ As in the C-W-B
method, Bourgain[Bo5] uses and develops the technique by
Fr\"ohlich and   Spencer[F-S] to investigate the inverse of $A$.
Even if one gets  the inverse $A^{-1}$, one has to show that
$A^{-1}$ is off-diagonal decay, in order to control the
exponent-weight norm of $A^{-1}$. This technique depends heavily
on the algebraic structure of $(k,\omega)+\Omega_j$. For nonlinear
wave equation of spatial $d$-dimension, $$\Omega^0_j\approx
\sqrt{|j|^2}=\sqrt{j_1^2+...+j_d^2}\ ,\ \ j=(j_1,...,j_d).$$ This
makes the algebraic structure of $(k,\omega)+ \Omega_j$ very
intricate, if $d>1$.

In the present paper, we consider Hamiltonian (1.2) with
$\iota=\infty$. Following Bourgain's idea in [Bo5], we put (1.7)
into the ``integrable" part (1.10) and find a symplectic
transformation to eliminate (1.6) and (1.5). This leads to solve
the homological equation (1.12) of variable coefficients. Thus, we
also have   to investigate the inverse of  the ``big" matrix $A$
as in [Bo5]. However, we do not use the technique by Fr\"ohlich
and Spencer[F-S]. We impose a symmetry condition on $R$ as follows
$$R(-x,y,z,\bar z)=R(x,y,\bar z,z)\tag 1.14$$ where $(x,y)$ is in some
subset of $\Bbb T^n\times\Bbb R^n$ and $z,\bar z$ are in some
Hilbert space. With the symmetry condition (1.14) we can prove
that the spectra of $A$ are ``twisted" with respect to $\omega$,
thus that the spectra are different zero if digging out some
$\omega$'s of small Lebesgue measure. This implies the existence
of $A^{-1}$. In addition, we do not need to know that $A^{-1}$ is
off-diagonal decay,  by choosing Sobolev space as our working
space instead of an exponential weight space. Generally it is very
hard to prove that $A^{-1}$ is off-diagonal decay, which involves
the intricate algebraic structure of $(k,\omega)+\Omega_j$.  It
should be noted that the symmetry condition (1.14) is fulfilled by
many partial differential equations, such as nonlinear wave
equation (1.1), nonlinear Schr\"odinger equation, nonlinear beam
equation and others.

The rest of the present paper is organized as follows. In $\S$2.,
a KAM theorem to deal Hamiltonian (1.2) with $\iota=\infty$ is
given out. This theorem does not require the second Melnikov's
conditions. In $\S$ 3., the KAM theorem is used to obtain
invariant tori and quasi-periodic solution for nonlinear wave
equation and nonlinear Schr\"odinger equation. The $\S$4-7 are
devoted to the proof of the KAM theorem. The fifth section  is the
essential part of the present paper. \vskip0.2cm \noindent{\bf
Acknowledgement. } {\it The present author is very grateful to
Professor D. Bambusi for encouragement.}

 \subhead 2. A  KAM theorem \endsubhead

{\it 2.1. Some notations.} Denote by $(\ell_2,||\cdot||)$ the
usual space of the square summable sequences, and  by
$(L^2,||\cdot||)$ the space of the square integrable functions. By
$|\cdot|$ the Euclidian  norm.  Let  $p\ge d/2$. For a sequence
$z=(z_j\in\Bbb C: j\in\zd)$, we define its norm as follows:
$$||z||_{p}^2=\sum_{j\in\zd}|j|^{2p}|z_j|^2.\tag 2.1$$
Let $\ell^{p}$ be the set of all sequences satisfying (2.1). It is
easy to see that $\ell^{p}$ is a Hilbert space with an  inner
product corresponding to (2.1).
 Introduce the  phase space:
$$\Cal P:=(\Bbb C^n/2\pi\Bbb Z^n)\times \Bbb
C^n\times\ell^{p}\times\ell^{p},\tag 2.2$$ where $n$ is a given
positive integer.
   We endow $\Cal P$
with a symplectic structure
$$dx\wedge dy+\sqrt{-1}dz\wedge d\bar z=dx\wedge dy+\sqrt{-1}\sum_{j\in\zd}dz_j\wedge d\bar z_j,\quad (x,y,z,\bar z)\in \Cal
P.\tag 2.3
$$
Let
$$\Cal T^n_0=(\Bbb R^n/2\pi\Bbb
Z^n)\times\{y=0\}\times\{z=0\}\times\{\bar z=0\}\subset\Cal P.\tag
2.4
$$ Then $\Cal T^n_0$ is a torus in $\Cal P$. Introduce a complex
neighborhoods of $\Cal T^n_0$ in $\Cal P$:
$$ D(s,r):=\{(x,y,u)\in  \Cal P:|\text{Im} x|<s, |y|<r^2,
||z||_{p} <r,||\bar z||_{p} <r\}\tag 2.5$$ where $r,s>0$ are
constants.

For $\tilde r>0$ we define the weighted phase norms
$$_{\tilde r}\W W\W_{p} =|X|+\frac{1}{\tilde r^2}|Y|+
\frac{1}{\tilde r}||Z||_{p}+\frac{1}{\tilde r}||\bar Z||_{p} \tag
2.6
$$
 for $W=(X,Y,Z, \bar Z)\in \Cal P$.
 Let $\Cal O\subset \Bbb R^ {n}$ be compact and of
positive Lebesgue measure. For a map $W:D(s,r)\times \Cal O\to
\Cal P$, set
$$_{\tilde r}\W W\W_{p,D(s,r)\times\Cal O}:=
\sup_{(x,\xi)\in D(s,r)\times\Cal O}\  _{\tilde r}\W
W(x,\xi)\W_{p} \tag 2.7$$ and
$$\ _{\tilde r}\W W\W_{ p,D(s,r)\times\Cal O}^{\Cal L}:=
\max_{1\le j\le n}\sup_{ D(s,r)\times\Cal O} \ _{\tilde r}\W
\p_{\xi_j}W(x,\xi)\W_{ p},\ \ \xi=(\xi_1,...,\xi_n). \tag 2.8$$

 Denote by $\Cal
L(\ell^{p},\ell^{ p})$ the set of all bounded linear operators
from $\ell^{p}$ to $\ell^{ p}$ and by $|||\cdot|||_{p}$ the
operator norm. For any subset $S\subset\zd$ with cardinality
$S^\sharp<\infty$ and a finitely dimensional vector $u=(u_j\in\Bbb
C:j\in S)$ and a matrix $U=(U_{ij}\in\Bbb C:i,j\in S)$ of finite
order, let
$$\tilde u=(\tilde u_j:j\in\zd),\quad \text{here }\tilde u_j=\cases & u_j, \ j\in S\\& 0,\ j\in\zd\setminus
S,\endcases$$ and
$$\tilde U=(\tilde U_{ij}:i,j\in\zd),\quad \text{here }
\tilde U_{ij}=\cases & U_{ij}, \ i, j\in S\\& 0,\ i \text{ or }
j\in\zd\setminus S.\endcases$$ Define $||u||_p:=||\tilde u||_p$
and $|||U|||_p:=|||\tilde U|||_p$. Similarly define
$||u||:=||\tilde u||$ and $|||U|||:=|||\tilde U|||$. The following
lemma states the relation between $|||\cdot|||_p$ and
$|||\cdot|||$.

\proclaim{Lemma 2.1} (i). Let $Y=(Y_{ij}:|i|,|j|\le K)$ be
$K\times K$ matrix with some $K>0$. Then $|||Y|||_p\le K^p|||Y|||$
and $|||Y|||\le |||Y|||_p$. (ii). Let $Y=(Y_j:|j|\le K)$ be a
$K$-dimensional vector. Then $||Y||_p\le K^p||Y||$ and $ ||Y||\le
||Y||_p$.
\endproclaim
\demo{Proof} The proofs are trivial. We omit it.\enddemo

 In the whole of this paper, by $C$ or $c$ a universal constant, whose
size may be different in different place. If $f\le C g$, we write
this inequality as $f\lessdot g$ when we dot not care the size of
the constant $C$. Similarly, if $f\ge C g$ we write $f\gtrdot g$.

\vskip 0.2cm

{\it 2.2. The statement of the KAM theorem.}  For two vectors
$b,c\in \Bbb C^k$ or $\Bbb R^k$, we write
$(b,c)=\sum_{j=1}^kb_jc_j$ if $k<\infty$. If $k=\infty$, we write
$\<b,c\>=\sum_{j=1}^\infty b_jc_j$. Consider an infinitely
dimensional Hamiltonian in the parameter dependent normal form
$$N_0=(\omega^0(\xi),y)+ \sum_{j\in \Bbb Z^d}\Omega_j^0(\xi)z_j\bar z_j,
\quad (x,y,z,\bar z)\in \Cal P.\tag 2.9$$
  The tangent
frequencies $\omega^0=(\omega_1^0,\cdots,\omega_n^0)$ and the
normal frequencies $\Omega_j^0$'s ($j\in \zd$) depend on $ n$
parameters $\xi\in \Cal O_0\subset \Bbb R^{n}$,where $\Cal O_0$ a
given compact set of positive Lebesgue measure. Let
$$\Omega^0(\xi)=\text{diag}(\Omega^0_j(\xi):j\in\zd).\tag 2.10$$
Let
$$N_0=(\omega^0(\xi),y)+ \<\Omega^0(\xi)z,\bar z\>,\tag 2.11$$ where
 $z=(z_j:j\in\zd)$.
% Let $$\Cal I=\pmatrix 0 & 1\\ -1 & 0 \endpmatrix,\ \ J=\sqrt{-1}\Cal I\tag 2.10$$ and
%$$J_m=\text{diag}(\underbrace{J,...,J}_{m}),\quad J_\jmath=\text{diag}(\underbrace{J,...,J}_{\jmath^\sharp}),\quad J_\infty=\text{diag}(\underbrace{...,J,...,J,...}_{\infty}).\tag 2.11$$
 The
Hamiltonian equation of motion of $N_0$ are
$$\dot x=\omega^0(\xi),\quad  \dot y=0,\quad \dot
z=\sqrt{-1}\Omega^0(\xi)z,\quad \dot {\bar
z}=-\sqrt{-1}\Omega^0(\xi)\bar z. \tag 2.12$$
 Hence, for each $\xi\in \Cal O_0$, there is an invariant
$n$-dimensional torus $$\Cal T^n_0=\Bbb T^n\times
\{y=0\}\times\{z=\bar z=0\}$$ with frequencies $\omega^0(\xi)$.
The aim is to prove the persistence of the torus $\Cal T^n_0$, for
``most" (in the sense of Lebesgue measure) parameter vector
$\xi\in \Cal O_0$, under small perturbation $R$ of the Hamiltonian
$N_0$. To this end the following assumptions are required.

\noindent{\it\bf Assumption A:} {\it (Multiplicity.)  Assume that
there are constants $c_1,c_2>0$ such that for all $\xi\in\Cal
O_0$,
$$(\Omega_j^0)^\sharp\le c_1|j|^{c_2}\tag
2.13$$ where we denote by $(\cdot)^\sharp$ the multiplicity of
$(\cdot)$.}

 \noindent{\it\bf Assumption B:} {\it (Non-degeneracy.)
There are two absolute constant $c_3,c_4>0$ such that
$$\sup_{\xi\in\Cal O_0}|\text{det }\p_{ \xi}\omega^0(\xi)|\ge c_3,\quad
\sup_{\xi\in\Cal O_0}|\p_{ \xi}^j\omega^0(\xi)|\le c_4,\
j=0,1.\tag 2.15$$

 \noindent{\it\bf Assumption C:} {\it (Analyticity of parameters.)} Assume that both $\omega^0(\xi)$ and
$\Omega^0(\xi)$ are analytic in each entry $\xi_l$  of the
variable vector $\xi=(\xi_1,...,\xi_l,...,\xi_{n})\in\Cal O_0$ and
assume that both $\omega^0(\xi)$ and $\Omega^0(\xi)$ are real for
real argument $\xi$.

\noindent{\it\bf Assumption D:} {\it (Bounded conditions of Normal
frequencies.) Assume that there exists constants $c_5,c_6,c_7
>0$ and constant $\kappa>0 $ such that
$$\inf_{\xi\in\Cal O_0}\Omega^0_j \ge c_5|j|^\kappa+c_6,\tag
2.16$$$$ \sup_{\xi\in\Cal O_0}|\p_{\xi_l}\Omega^0_j|\le c_7\ll 1,\
l=1,...,n \tag 2.17$$ uniformly for all $j$.}

\noindent{\it\bf Assumption E:} {\it (Regularity.) Let $s_0,r_0$
given. Assume the perturbation term $ R^0(x,y,z,\bar z;\xi)$ which
is defined on the domain $D(s_0,r_0)\times\Cal O_0$ is analytic in
the space coordinates and also analytic in each entry $\xi_l$
($l=1,..., n$) of the parameter vector $\xi\in\Cal O_0$, and is
real for real argument, as well as, for each $\xi\in \Cal O_0$
 its Hamiltonian vector field $X_{ R^{0}}:=(
R^{0}_y,- R^{0}_x, \sqrt{-1}\p_z R^{0},-\sqrt{-1}\p_{\bar z}
R^{0})^T$ defines
 a analytic map
$$X_{ R^{0}}:D(s_0,r_0)\subset\Cal P\to\Cal P
$$ where  $T\equiv$
transpose. Also assume  that $X_{ R^{0}}$ is analytic in  each
entry of  $\xi\in\Cal O_0$.

\noindent{\it\bf Assumption F:} {\it (Symmetry.) For any
$(x,y,z,\bar z,\xi)\in D(s_0,r_0)\times\Cal O_0$,
$$R^0(-x,y,z,\bar z,\xi)=R^0(x,y,\bar z, z,\xi).\tag 2.18$$

\proclaim{Theorem 2.1} Suppose $H=N_0+ R^0$ satisfies assumptions
A--F and smallness assumption:
$$\ _{r_0}\W X_{
R^{0}}\W_{ p,D(s_0,r_0)\times\Cal O_0} <\e,\quad
  _{r_0}\W X_{
R^{0}}\W_{ p,D(s_0,r_0)\times\Cal O_0}^{\Cal L}< \e^{1/3}.\tag
2.19$$ Then, for  $0<\epsilon\ll 1$, there is a subset $\Cal
O_\e\subset \Cal O_0$ with $$\text{ Meas } \Cal O_\e\ge (\text{
Meas }\Cal O_0)(1-O(\epsilon)),$$ and there are a family of torus
embedding $\Phi:\Bbb T^n\times\Cal O_\e\to \Cal P$ and a map
$\omega_*:\Cal O_\e\to\Bbb R^n$ where $\Phi(\cdot,\xi)$ and
$\omega_*(\xi)$ is analytic in each entry $\xi_j$ of the parameter
vector $\xi=(\xi_1,...,\xi_n)$ for other arguments fixed, such
that for each $\xi\in \Cal O_\e$ the map $\Phi$ restricted to
$\Bbb T^n\times\{\xi\}$ is a analytic embedding of a rational
torus with frequencies $\omega_*(\xi)$ for the Hamiltonian $H$ at
$\xi$.

Each embedding is real analytic  on $\Bbb T^n\times\{\xi\}$, and
$$\ _{r_0}\W \Phi-\Phi_0\W_{ p,\Bbb T^n\times\Cal O_\e}\le c\e,\quad
\ _{r_0}\W\Phi-\Phi_0\W_{p, \Bbb T^n\times\Cal O_\e}^{\Cal L}\le
c\epsilon^{1/3},$$
$$|\omega_*-\omega|_{\Cal O_\e}\le c\e,\quad |\omega_*-\omega|_{\Cal O_\e}^{\Cal L}
\le c\epsilon^{1/3},$$
 where $\Phi_0$ is the
trivial embedding $\Bbb T^n \times\Cal O_0\to \Bbb
T^n\times\{y=0\}\times \{z=\bar z=0\}$, and $c>0$ is a constant
depending on $n$.
\endproclaim

\subhead 3. Application to nonlinear partial differential
equations
\endsubhead

 \noindent{\it 3.1 Application to nonlinear wave equations.}
Consider nonlinear wave equation
$$u_{tt}-\triangle u+M_\sigma u+\e u^3=0,\ \ \theta\in\Bbb T^d,\ d\ge 1\tag 3.1$$
where $u=u(t,\theta)$ and $\triangle=\sum_{j=1}^d\p_{\theta_j}^2$
and $M_\sigma$ is a real Fourier multiplier
$$M_\sigma\cos (j,\theta)=\sigma_j\cos(j,\theta),\
M_\sigma\sin (j,\theta)=\sigma_j\sin(j,\theta),\ \ \sigma_j\in\Bbb
R,\ j\in\zd.\tag 3.2$$ Pick a set $\un
e=\{e_1,...,e_n\}\subset\zd$. Let
$\un\zd=\zd\setminus\{e_1,...,e_n\}$. Following Bourgain[Bo3], we
assume
$$\cases &\sigma_{e_j}=\sigma_j,\ \ (j=1,...,n)\\ & \sigma_j=0,\ \ j\in
\un\zd.\endcases\tag 3.3$$  Let $\lambda_j^\pm$ ($j\in\zd$) be the
eigenvalues of the self-adjoint operator $-\triangle +M_\sigma$
subject to periodic b. c. $\theta\in\Bbb T^d$, and let
$\phi_j^\pm=\phi_j^\pm(\theta)$ be the normalized eigenfunctions
corresponding to $\lambda_j^\pm$. Then
$$\lambda_j^\pm=|j|^2+\sigma_j=\sum_{l=1}^dj_l^2+\sigma_j,\
j=(j_1,...,j_d)\in\zd,\tag 3.4$$ and
$$\phi_j^+=\frac{\sqrt 2}{(\sqrt{2\pi})^d}\cos(j,\theta),\
\phi_j^-=\frac{\sqrt 2}{(\sqrt{2\pi})^d}\sin(j,\theta).\tag 3.5$$
Note that $\{\phi_j^\pm:j\in\zd\}$ is complete orthogonal system
in $L^2(\Bbb T^d)$. Let
$$u(t,\theta)=\sum_{j\in\zd}q_j^\pm(t)\phi_j^\pm(\theta).\tag 3.6$$
Inserting (3.6) into (3.1), then
$$\ddot q_j^\pm+\lambda_j^\pm q_j^\pm+\e (u^3,\phi_j^\pm)_{L^2}=0,\
j\in\zd.\tag 3.7$$ Let $ \un q_j^\pm=\sqrt{\lambda_j^\pm}\dot
q_j^\pm$. Then
$$\cases &\dot q_j^\pm=\sqrt{\lambda_j^\pm}\un q_j^\pm\\&
\dot {\un q}_j^\pm=\frac{1}{\sqrt{\lambda_j^\pm}}\ddot
q_j^\pm=-\sqrt{\lambda_j^\pm}q_j^\pm-\frac{\e}{\sqrt{\lambda_j^\pm}}
(u^3,\phi_j^\pm).\endcases\tag 3.8$$ For simplifying notation, we
do not distinguish $+$ sign and $-$ sign in (3.8). For example, we
write $q_j^\pm$ as $q_j$. However, we should keep in mind that
$q_j$ runs over the set $\{q_j^+,q_j^-\}$. The system (3.8) is a
Hamiltonian system with its Hamiltonian function
$$H=H(q,\un q)=\frac12\sum_{j\in\zd}\sqrt{\lambda_j}(\un q_j^2+q_j^2)
+\e G(q),\tag 3.9$$ where
$$G(q)=\sum_{i,j,k,l\in\zd}G_{ijkl}q_iq_jq_kq_l\tag 3.10$$ and
$$G_{ijkl}=\frac{1}{\sqrt{\lambda_i\lambda_j\lambda_k\lambda_l}}
\int_{\Bbb T^d}\phi_i\phi_j\phi_k\phi_l dx.\tag 3.11$$ Let
$$\omega_j^0=\sqrt{\lambda_{e_j}}=\sqrt{|e_j|^2+\sigma_j},\ \
(j=1,...,n),\tag 3.12$$ and
$$\Omega_j^0=\sqrt{\lambda_j}=\sqrt{|j|^2},\ \ j\in
\un{\zd}.\tag 3.13$$ Introduce a symplectic coordinate change
$$\cases & p_{e_j}=\sqrt{2(1+y_j)}\sin x_j,\ \
q_{e_j}=\sqrt{2(1+y_j)}\cos x_j,\ \ j=1,...,n,\\ &
q_j=\frac{z_j+\bar z_j}{\sqrt 2},\quad p_j=\frac{z_j-\bar
z_j}{\sqrt{-1}\sqrt 2},\ j\in\un\zd.\endcases\tag 3.14$$ To
shorten notation, let $$\tilde y_i=\sqrt{2(1+y_i)}.$$ Then (3.9)
is changed into
$$H=H(x,y,z,\bar z)=(\omega,y)+\sum_{j\in\un\zd}\Omega_j^0z_j\bar z_j
+\e R^0\tag 3.15$$ where
$$\aligned R^0=&R^0(x,y,z,\bar z)\\=& \sum_{i,j,k,l\in\un e}G_{ijkl}\tilde y_i\tilde y_j
\tilde y_k\tilde y_l\cos x_i\cos x_j\cos x_k\cos x_l\\&+
 \sum_{i,j,k\in\un e, l\in\un\zd}G_{ijkl}\tilde y_i\tilde y_j\tilde y_k
 \cos x_i\cos x_j\cos x_k\frac{z_l+\bar z_l}{\sqrt 2}\\&
 +\sum_{i,j\in\un e, k,l\in\un\zd}G_{ijkl}\tilde y_i\tilde y_j
  \cos x_i\cos x_j\frac{z_k+\bar z_k}{\sqrt 2}
   \frac{z_l+\bar z_l}{\sqrt 2}\\&
  +\sum_{i\in\un e, j,k,l\in\un\zd}G_{ijkl}\tilde y_i
  \cos x_i\frac{z_j+\bar z_j}{\sqrt 2}\frac{z_k+\bar z_k}{\sqrt 2}
   \frac{z_l+\bar z_l}{\sqrt 2}\\&
   +\sum_{i, j,k,l\in\un\zd}G_{ijkl}\frac{z_i+\bar z_i}{\sqrt 2}\frac{z_j+\bar z_j}{\sqrt 2}\frac{z_k+\bar z_k}{\sqrt 2}
   \frac{z_l+\bar z_l}{\sqrt 2} \endaligned\tag 3.16$$
It follows that in a suitable domain of $(x,y,z,\bar z)$,
$$R^0(-x,y,z,\bar z)=R^0(x,y,\bar z, z).\tag 3.17$$ This implies
the symmetry assumption F is fulfilled.  By (3.5) and (3.11), it
is not difficult to verify that
$$G_{ijkl}=0,\ \ \text{ unless } i\pm j\pm k\pm l=0, \tag 3.18$$
for some combination of plus and minus signs. For $w,z\in\ell^p$,
the convolution $z\star w$ is defined by $(z\star
w)_j=\sum_{k\in\zd}w_{j-k}z_k$.  \proclaim{ Lemma 3.1} If $p>d/2$,
then $||w\star z||_p\le c ||w||_p||z||_p$ for $w,z\in\ell^p$ with
a constant $c$ depending only on $p$.
\endproclaim
\demo{Proof} The proof is very elementary. See
[p.294-295,P2].\enddemo Let $s_0=r_0=1$,
$\omega^0=(\omega_1,...,\omega_n)$ and
$\Omega^0=(\Omega_j^0:j\in\un\zd)$. And let parameter
$\sigma=(\sigma_1,...,\sigma_n)$ runs over $\Cal O_0:=[1,2]^d$.
Note that $R^0$ is independent of $\sigma$. By (3.18) and Lemma
3.1, we see that $R^0$ is well defined on $D(s_0,r_0)$ and  the
Assumption E holds true; moreover,$$ _{r_0}\W X_{\e
R^0}\W_{p,D(s_0,r_0)\times\Cal O}\lessdot \e,\ _{r_0}\W X_{\e
R^0}\W_{p,D(s_0,r_0)\times\Cal O}^\Cal L=0.\tag 3.19$$ This
implies that (2.19) is fulfilled. It follows from (3.12) that the
Assumption B (that is, (2.15) ) holds true. By (3.13), we see that
the assumptions A, C, D are fulfilled. In particular, (2.16) holds
true with $\kappa=1$. Using Theorem 2.1, we have the following
theorem.

\proclaim{Theorem 3.1} There is a subset $\un{\Cal O}_0\subset\Cal
O_0$ with $\text{Meas }\un{\Cal O}\ge (1-C\e)\text{Meas }{\Cal O}$
such that for any $\sigma\in\un{\Cal O}_0$, Eq.(3.1) with small
$\e$ has a rotational invariant torus of frequency vector
$\omega^0=\omega^0(\sigma)$. The motions on the torus are
quasi-periodic with frequency $\omega^0$.
\endproclaim

\noindent{\it 3.2 Application to nonlinear Schr\"odinger
equations.} Consider nonlinear Schr\"odinger equation
$$\sqrt{-1}u_{t}-\triangle u+M_\sigma u+\e u|u|^2=0,\ \
\theta\in\Bbb T^d,\ d\ge 1,$$ where  $M_\sigma$ is a Fourier
multiplier
$$M_\sigma e^{\sqrt{-1}(j,\theta)}=\sigma_je^{\sqrt{-1}(j,\theta)},\ \ \sigma_j\in\Bbb
R,\ j\in\zd.$$ \proclaim{Theorem 3.2} Assume $\sigma_j$'s satisfy
(3.3).  There is a subset $\un{\Cal O}_0\subset\Cal O_0$ with
$$\text{Measure }\un{\Cal O}\ge (1-C\e)\text{Measure }{\Cal O}$$
such that for any $\sigma=(\sigma_1,...,\sigma_n)\in\un{\Cal
O}_0$, the nonlinear Schr\"odinger equation with small $\e$ has a
rotational invariant torus of frequency vector
$\omega^0=\omega^0(\sigma)$ where
$$\omega_j^0=|e_j|^2+\sigma_j,\ j=1,...,n.$$ The motions on the torus are
quasi-periodic with frequency $\omega^0$.
\endproclaim
\demo{Proof} The proof is similar to that of Theorem 3.1. Note
that in this case, (2.16) is fulfilled for $\kappa=2$. We omit the
details.
\enddemo
 %%%%%%%%%%%%
%%%%%%%%%%%%
%%%%%%%%%%%%
%%%%%%%%%%%%
%%%%%%%%%%%%
%%%%%%%%%%%%
%%%%%%%%%%%%
 \subhead 4. The linearized equation\endsubhead

{\it 4.1. Unperturbed linear system. } Recall $\Omega_j^0$'s
satisfy the {\bf Assumptions A, C, D}.  Assume that there is a new
tangent frequency vector   \footnote{ Here $\omega=\omega(\xi)$ is
not necessary to be $\omega^0$, due to the well-known phenomenon
the shift of frequency in the  KAM iteration.} $\omega$ satisfing
the {\bf Assumptions B, C}. Let $\Cal O\subset\Cal O_0$ be an open
set. In this section, we pick two large constants $K_-$ and $K$.
And let $s=1/K_-$ and $s^\prime=4/K$. (In the $m$-th KAM step, we
will choose $K_-=K_{m-1}\approx 2^{(m-1)^2/2}$ and $K=K_{m}\approx
2^{m^2/2}$.) Let
$$D(s)=\{x\in\Bbb C^n/2\pi\Bbb Z^n:|\Im x|<s \},\ \ s=1/K_-.$$
 Let $X_{\un N}$ be a linear  Hamiltonian system with Hamiltonian function:
$$\aligned \un N=&(\un\omega(\xi),y)+\sum_{j\in\zd}\Omega_j^0(\xi)z_j\bar z_j
+\\& +\frac12\sum_{i,j\in\zd}\un B_{ij}^{zz}(x;\xi)z_i
z_j+\sum_{i,j\in\zd}\un B_{ij}^{z\bar z}(x;\xi)z_i\bar
z_j+\frac12\sum_{i,j\in\zd}\un B_{ij}^{\bar z\bar z}(x;\xi)\bar
z_i\bar z_j)\\:=&(\un\omega,y)+(\Omega^0z,\bar z)+\frac12\<\un
B^{zz}z, z\>+\<\un B^{z\bar z}z, \bar z\>+\frac12\<\un B^{\bar
z\bar z}\bar z,\bar z\>\endaligned \tag 4.1$$ where  $\un
B_{ij}(x;\xi)$'s are analytic in $x\in D(s)$ for any fixed
$\xi\in\Cal O$, and all of $\omega(\xi),\Omega^0_\jmath(\xi)$ and
$ B_{ij}(x;\xi)$ (for fixed $x\in D(s)$) are analytic in each
entry $\xi_l$ ($l=1,..., n$) of $\xi\in\Cal O$. Assume $B_{ij}$'s
satisfy the following conditions:

\par\noindent {\it (1.) Symmetry.}
$$\un B^{zz}_{ij}(x,\xi)=\un B^{zz}_{ji}(x,\xi),\ \un B^{zz}_{ij}(-x,\xi)=
\un B^{\bar z\bar z}_{ij}(x,\xi) \tag 4.2 $$
$$\un B^{z\bar z}_{ij}(x,\xi)=\un B^{z\bar z}_{ji}(x,\xi),\
\un B^{z\bar z}_{ij}(-x,\xi)=\un B^{ z\bar z}_{ij}(x,\xi)\tag 4.3
$$ for any $(x,\xi)\in D(s)\times\Cal O$.
\par\noindent {\it (2.) Finiteness of Fourier modes.}
$$\un B_{ij}^{zz}(x;\xi)=\sum_{|k|\le K}\widehat{\un B_{ij}^{zz}}(k)
e^{\sqrt{-1}(k,x)},\cdots \tag 4.4$$ where $K\gg 1$ is a constant
independent of $x,\xi,i$ and $j$. (In the $m$-th KAM iteration, we
will take $K\approx 2^{m^2/2}$.)
\par\noindent {\it (3.) Boundedness.}
$$\sup_{D(s)\times\Cal O}|||\un B(x;\xi)|||^*_p\lessdot  \e\ll 1,\tag 4.5 $$
where $*\equiv$ the blank and $\Cal L$.
\par\noindent {\it (4.) Non-degenerateness.}
 Moreover, we assume $$c_1\ge\left|\text {det}\frac{\p
\omega(\xi)}{\p\xi}\right|\ge c_2>0,\  \forall \ \xi\in \Cal
O.\tag 4.6$$ (This assumption enables  us to regard $\omega$
instead of $\xi$ as a parameter vector.)

\par{\it 2. split and estimate for small perturbation.} We
now consider a perturbation $$H=\un N+\grave{R}$$ where
 $\grave{R}=\grave{R}(x,y,u;\xi)$
 is a Hamiltonian defined on
 $D(s,r)$ and  depends on the parameter $\xi\in\Cal O$. Assume that the vector field
$X_{\grave R}:\quad D(s,r)\times\Cal O\to\Cal P$ is real for real
argument, analytic in $(x,y,z,\bar z )\in D(s,r)$ and
  in  each entry of $\xi\in \Cal O$.
 We
assume that there are quantities\footnote{We will take
$\varepsilon=\e^{(4/3)^m}$ and $\varepsilon^{\Cal
L}=\varepsilon^{1/3}$ in the $m$-th KAM iteration step.}
$\varepsilon=\varepsilon(r,s,\Cal O)$ and $\varepsilon^{\Cal
L}=\varepsilon^{\Cal L}(r,s,\Cal O)$ which depend on $r,s,\Cal O$
 such that
$$\ _{r}\W X_{\grave R}\W_{p,D(s,r)\times\Cal O}\lessdot \varepsilon,\quad
 \ _r\W X_{\grave R}\W_{p,D(s,r)\times\Cal O}^{\Cal L}\lessdot
 \varepsilon^{\Cal L},\ \ \varepsilon<\varepsilon^{\Cal L}\ll 1.\tag 4.7$$
 Let
$$R=\sum\Sb2|m|+|q_1+ q_2|\le
2\\m\in\zd, q_1,q_2\in\Bbb Z^\infty
\endSb\sum_{k\in\zd,}\widehat{\grave R_{mq_1q_2}}(k)
e^{\sqrt{-1}(k,x)}y^mz^{q_1} \bar z^{ q_2},\tag 4.8$$ with the
Taylor-Fourier coefficients $\widehat {\grave R_{mq_1q_2}}(k) $ of
$\grave{R}$ depending on $\xi\in\Cal O$, and being analytic in
each entry $\xi_j$ of $\xi$.  We see that $R$ is a partial
Taylor-Fourier expansion of $\grave R$. We will approximate
$\Grave {R}$ by $R$.   Now we give some estimates of $R$.
\proclaim{Lemma 4.1}
$$\ _r\W X_{R}\W_{p,D(s,r)\times
\Cal O}^*\lessdot \ _r\W X_{\grave{R}}\W_{p,D(s,r)\times\Cal O}^*
\le \varepsilon^*,\tag 4.9
$$
$$\ _{\eta r}
\W X_{R}-X_{\grave R}\W_{p,D(s,4\eta r)\times\Cal O}^*
 \lessdot \eta \cdot \ _r\W X_{\grave
R}\W_{p,D(s,r)\times\Cal O}^* \lessdot \eta \varepsilon^*, \tag
4.10$$ for any $ 0<\eta\ll 1$, where $^*=$ the blank or $\Cal L$,
for example, $\varepsilon^*=\varepsilon$ or $\varepsilon^{\Cal
L}$.
\endproclaim
\demo{Proof} The proof is similar to that of  formula (7) of
[P1,129].
\enddemo
With this lemma, we decompose $R=R^0+R^1+R^2$, where
$$\aligned & R^0=R^x+(R^y,y),\\
& R^1=\<R^{z},z\>+\<R^{\bar z},\bar z\>,\\
& R^2=\frac12(\<R^{zz}z,z\>+\<R^{z\bar z}z,\bar z\>+\<R^{\bar
z\bar z}\bar z,\bar z\>),\endaligned$$ with $R^x,R^y:
D(s)\times\Cal O\to \Bbb C^n$; $ R^{z}, R^{\bar z}:D(s)\times\Cal
O\to \ell^p$; $R^{zz}, R^{z\bar z}, \text{ etc. }D(s)\times\Cal
O\to \Cal L(\ell^p,\ell^p)$. For any vector or matrix
$$Y=\sum_{k,m\in\zd;q,\bar q\in\Bbb Z^\infty}\widehat{Y_{m,q,\bar q}}(k)e^{\sqrt{-1}(k,x)}y^mz^q\bar
z^{\bar q}$$ we introduce the cut-off operator $\Gamma_K$ as
follows:
$$\Gamma_KY=\sum\Sb k,m\in\zd;q,\bar q\in\Bbb Z^\infty\\ |k|\le K\endSb\widehat{Y_{m,q,\bar q}}(k)e^{\sqrt{-1}(k,x)}y^mz^q\bar
z^{\bar q}\tag 4.11$$

\proclaim{Lemma 4.2} Assume $(s-s^\prime)K\ge |\ln\eta|$. Let
$R_K=R-\Gamma_KR$. We have
$$\ _{r}\W X_{(\Gamma_KR)}\W_{p,D(s^\prime,r)\times\Cal O}^*\le \ _r\W X_{R}\W_{p,D(s,r)\times
\Cal O}^*\lessdot \varepsilon^*,\tag 4.12$$ $$ \ _{r}\W
X_{R_K}\W_{p,D(s^\prime,r)\times\Cal O}^{*}
\lessdot\eta\varepsilon^{*}\tag 4.13$$ where  $^*=$ the blank or
$\Cal L$.
\endproclaim
\demo{Proof} The proof of (4.12) is obvious. Let us give the proof
of (4.13). Write
$$R_K=R_K^x+(R_K^y,y)+\<R_K^z,z\>+\<R_K^{\bar z},\bar z\>+\<R_K^{zz}z,z\>+
\<R_K^{z\bar z} z,\bar z\>+\<R_K^{\bar z\bar z}\bar z,\bar z\>.$$
Note that the terms $R_K^x$, $R_K^y$, and so on, are analytic in
$x\in D(s)$ for fixed $\xi\in\Cal O$. And observe that
$$R_K^x=\sum_{|k|>K}\widehat{R_K^x}(k)e^{\sqrt{-1}(k,x)},\cdots. $$
 Then by Cauchy's formula,
we have $|\widehat{R_K^x}(k)|\le e^{-s|k|}\sup_{D(s)}|R^x|$, and
so on. Thus,
$$\ _{r}\W X_{R_K}\W_{p,D(s^\prime,r)\times\Cal O}^*\le \ _r\W
X_R\W_{p,D(s,r)\times\Cal O}^*\sum_{|k|>K}e^{-\sigma |k|}\le
\eta\varepsilon^*.$$  $\qed$
\enddemo
Now we can write
$$\grave R=\Gamma_KR+R_K+(\grave R-R).\tag 4.14$$
 Let $\Cal R=\Gamma_K R$. Hence we can write
 $$\aligned \Cal R=&\Cal R^x+(\Cal R^y,y)+\<\Cal R^z,z\>+
 \<\Cal R^{\bar z},\bar z\>\\&+\frac12(\<\Cal R^{zz}z,z\>+
 \<\Cal R^{z\bar z}z,\bar z\>+
 \<\Cal R^{\bar z\bar z}\bar z,\bar z\>).\endaligned\tag 4.15
  $$ Noting (4.8) and (4.11) we have
$$\Cal R^*=\sum_{|k|\le K}\widehat{R^*}(k)e^{\sqrt{-1}(k,x)},\tag 4.16 $$
where $*=x,y, z,\bar z, zz,z\bar z,\bar z\bar z$.

 \proclaim{Lemma
4.3} Under the smallness assumption (4.7) on $\grave R$, the
following estimates hold true:
$$|\p_x\Cal R^x|_{D(s)\times\Cal O}\le Kr^2\varepsilon,\quad
|\p_x\Cal R^x|^{\Cal L}_{D(s)\times\Cal O}\le Kr^2
\varepsilon^{\Cal L}\tag 4.17
$$
$$|\Cal R^y|_{D(s)\times\Cal O}\le K\varepsilon,\quad
|\Cal R^y|_{D(s)\times\Cal O}^{\Cal L} \le K\varepsilon^{\Cal L}
\tag 4.18$$
$$||\Cal R^u||_{ p, D(s)\times\Cal O}\le K r \varepsilon,\quad
||\Cal R^u||_{p, D(s)\times\Cal O}^{\Cal L} \le
Kr\varepsilon^{\Cal L},\ \ u\in\{ z,\bar z\} \tag 4.19$$
$$|||\Cal R^{uv}|||_{p,D(s)\times\Cal O} \le K\varepsilon
,\quad |||\Cal R^{uv}|||_{p,D(s)\times\Cal O}^{\Cal L} \le
K\varepsilon^{\Cal L},\ \ u,v\in\{z , \bar z \}.\tag 4.20
$$
 \endproclaim
\demo{Proof} Consider $\Cal R^{z\bar z}$. Observe that $R^{z\bar z
}=\p_z\p_{\bar z} R|_{z=\bar z=0}$. By the generalized Cauchy
inequality (See Lemma A.3 in [P1]),
$$\aligned & |||\Cal R^{z\bar z}|||_{p,D(s)\times\Cal O}\le K|||R^{z\bar z}|||_{p,D(s)\times\Cal O}
\le \frac{K}{r}||\p_{\bar z}R||_{p,D(s,r)\times\Cal O}\\& \le \
_r\W X_R\W_{p,D(s,r)\times\Cal O}K<K\varepsilon.\endaligned\tag
4.21$$ The remaining proof is simple. We omit the details. $\qed$
\enddemo

\proclaim{Lemma 4.4} Assume $\grave R$ satisfies the following
symmetric condition:
$$\grave R(-x,y,z,\bar z;\xi)= \grave R(x,y,\bar z,z;\xi),\ \ \forall\
(x,y,z,\bar z;\xi)\in D(s,r)\times \Cal O.\tag 4.22$$ Let
$R^{zz}_{ij},R^{z\bar z}_{ij},R^{\bar z\bar z}_{ij}$ be the
elements of the matrices $R^{zz},R^{z\bar z},R^{\bar z\bar z}$,
respectively. Then we have
$$R^x(-x)=R^x(x),\ \ R^y(-x)=R^y(x) \tag 4.23$$
$$R^z(-x)=R^{\bar z}(x),\ \  R^{\bar z}(-x)=R^{ z}(x) \tag 4.24$$
$$R^{zz}_{ij}(x)=R^{zz}_{ji}(x),\ R^{zz}_{ij}(-x)=R^{\bar z\bar z}_{ij}(x)
\tag 4.25-1 $$
$$R^{z\bar z}_{ij}(x)=R^{z\bar z}_{ji}(x),\
R^{z\bar z}_{ij}(-x)=R^{ z\bar z}_{ij}(x)\tag 4.25-2 $$
\endproclaim
\demo{Proof} Noting that
$$R^{zz}_{ij}(x)=\p_{z_i}\p_{z_j}\grave R(x,y,z,\bar z;\xi)|_{y=0,z=\bar z=0},$$
and $\p_{z_i}\p_{z_j}=\p_{z_j}\p_{z_i}$ we get the first equation
of (4.25-1). Applying  $\p_{z_i}\p_{z_j}$ to both sides of (4.23),
we get
$$\aligned &\p_{z_i}\p_{z_j}\grave R(-x,y,z,\bar z;\xi)|_{y=0,z=\bar
z=0}\\=&
 \p_{z_i}\p_{z_j}\grave R(x,y,\bar z, z;\xi)|_{y=0,z=\bar z=0}\\=&
 \p_{\bar z_i}\p_{\bar z_j}\grave R(x,y, z,\bar z;\xi)|_{y=0,z=\bar z=0},\endaligned$$
that is, the second equation of (4.25-1) holds. The remaining
proofs are similar to the previous one.$\qed$
\enddemo
Recall $\Cal R=\Gamma_K R$. Let
$$B^{zz}=\un B^{zz}+\Cal R^{zz}, B^{zz}=\un B^{zz}+\Cal R^{zz},
B^{\bar z\bar z}=\un B^{\bar z\bar z}+\Cal R^{\bar z\bar z}\tag
4.26$$
$$\omega=\un\omega+\widehat{R^y}(0).\tag 4.27$$
 and
$$N=(\omega,y)+(\Omega^0z,\bar z)+\frac12\< B^{zz}z, z\>+\<
B^{z\bar z}z, \bar z\>+\frac12\< B^{\bar z\bar z}\bar z,\bar z\>,
\tag 4.28$$ where $B^{z\bar z}=B^{\bar z z}$. Then by (4.25),
(4.2) and (4.3),
$$ \aligned & B^{zz}_{ij}(x,\xi)=B^{zz}_{ji}(x,\xi),\  B^{zz}_{ij}(-x,\xi)=
 B^{\bar z\bar z}_{ij}(x,\xi)  \\&
B^{z\bar z}_{ij}(x,\xi)= B^{z\bar z}_{ji}(x,\xi),\
 B^{z\bar z}_{ij}(-x,\xi)= B^{ z\bar z}_{ij}(x,\xi)\endaligned
 \tag 4.29
$$ for any $(x,\xi)\in D(s)\times\Cal O$.
Assume $K\varepsilon<\e$ and $K\varepsilon^\Cal L<\e$. By (4.5)
and (4.20), we have
$$|||B^{zz}|||_p^*, |||B^{z\bar z}|||_p^*, |||B^{\bar z\bar z}|||_p^*
\lessdot\e\tag 4.30$$ where $*=$ the blank or $\Cal L$.
 By (4.18) and $K\varepsilon\ll 1$ as
well as (4.6), we get that (4.6) still hold true for $\omega$:
$$\tilde c_1\ge\left|\text {det}\frac{\p
\omega(\xi)}{\p\xi}\right|\ge \tilde c_2>0,\  \forall \ \xi\in
\Cal O.\tag 4.31$$
 In view
of (4.27,28,14,15),
$$\aligned H=& N+\Cal R^x+(\Cal R^y-\widehat{R^y}(0),y)+\<\Cal R^z,z\>+\<\Cal R^{\bar z},\bar
z\>\\&+(R-\Gamma_KR)+(\grave R-R).\endaligned\tag 4.32$$

%%%%%%%%%%%%%%%%
%%%%%%%%%%%%%%%%
%%%%%%%%%%%%%%%%
%%%%%%%%%%%%%%%%

{\it 4.3. Derivation of homological equations.} The KAM theorem is
proven by the usual Newton-type iteration procedure which involves
an infinite sequence of symplectic coordinate changes. Each
coordinate change is obtained as the time-$1$ map $X_F^t|_{t=1}$
of a Hamiltonian vector field $X_F$. Its generating Hamiltonian
$F$ solves the linearized equation
$$\{F,N\}=\Cal R^x+(\Cal R^y-\widehat{R^y}(0),y)+\<\Cal R^z,z\>+\<\Cal R^{\bar z},\bar
z\>\tag 4.33$$ where $\{\cdot,\cdot\}$ is Poisson bracket with
respect to the symplectic structure $dx\wedge dy+\sqrt{-1}dz\wedge
d \bar z$. Without loss of generality, we assume $\widehat{\Cal
R^x}(0)=0$ since it does not affect the dynamics.
 We are now in position to find a solution of
the equation (4.33) and to give some estimates of the solution. To
this end, we suppose that $F$ is of the same form as the right
hand of (4.33), that is, $F=F^0+F^1$, where
$$\left\{\aligned & F^0=F^x+(F^y,y),\\
& F^1=\<F^{z},z\>+\<F^{\bar z},\bar z\>,\endaligned\right.\tag
4.34
$$ with $F^x,F^y,F^{z}, F^{\bar z}$ depending on $x,\xi$, and
$\widehat{F^x}(0)=0$, $\widehat{F^y}(0)=0$. Let
$\p_\omega=(\omega,\p_x)$ where $(\cdot,\cdot)$ is the usual inner
product in $\Bbb R^n$.
%%%%%%%%%%%
%%%%%%%%%%%
% 2006年5月10日0：10分修改到此
%%%%%%%%%%%
%%%%%%%%%%%
Using (4.34) and (4.28) we can compute the Poisson Bracket
$\{F,N\}$:
$$\aligned \{F,N\}=&-\frac12\<(F^y,\p_x)B^{zz}z, z\>-\<(F^y,\p_x)B^{z\bar z}
z,\bar z\> -\frac12\<(F^y,\p_x)B^{\bar z\bar z}\bar z,\bar
z\>\\&+\pw F^x+\<\pw F^y, y\>+\<\pw F^z,z\>+\<\pw F^{\bar z},\bar
z\>\\&+\sqrt{-1}(\<B^{zz}F^{\bar z},z\>+\<B^{z\bar z}F^{\bar
z},\bar z\>+\<\Omega^0F^z,z\>)\\&-\sqrt{-1}( \<B^{\bar z\bar z}F^{
z},\bar z\>+\<B^{z\bar z}F^{z},z\>+\<\Omega^0F^{\bar z},\bar z\>).
\endaligned\tag 4.35$$
 From (4.33,35) we derive the
homological equations:
$$\pw F^x   =\Gamma_K R^x(x,\xi),\tag 4.36$$
$$\pw F^y =\Gamma_K R^y(x,\xi)-\widehat{ R^y}(0),\tag 4.37$$
$$-\sqrt{-1}\pw F^z+\Gamma_K(\Omega^0F^z-B^{z\bar z}F^z+B^{zz}F^{\bar z})
=-\sqrt{-1}\Gamma_K R^z(x,\xi),\tag 4.38-1$$
$$\sqrt{-1}\pw F^{\bar z}+\Gamma_K(\Omega^0F^{\bar z}-B^{z\bar z}F^{\bar z}
+B^{\bar z\bar z}F^{ z})=\sqrt{-1}\Gamma_K R^{\bar z}(x,\xi),\tag
4.38-2$$
 We should note that if $F$ solves the
equations (4.36-38), then $F$ solves the following equation
$$\aligned \{F,N\}=&\frac12\<(F^y,\p_xB^{zz})z, z\>+\<(F^y,
\p_xB^{z\bar z}) z,\bar z\> +\frac12\<(F^y,\p_xB^{\bar z\bar
z})\bar z,\bar z\>\\&+ \sqrt{-1}\left\<(1-\Gamma_K)(B^{\bar z\bar
z}F^z),\bar z\right\>+ \sqrt{-1}\left\<(1-\Gamma_K)(B^{z\bar z}F^{
z}), z\right\>\\&-\sqrt{-1}\left\<(1-\Gamma_K)(B^{ z z}F^{\bar
z}), z\right\>- \sqrt{-1}\left\<(1-\Gamma_K)(B^{z\bar z}F^{\bar
z}), \bar z\right\>
\endaligned\tag 4.39 $$ instead of (4.35).

%%%%%%%%%%%%%%%%
 \proclaim{Lemma 4.5} If $F$ in some sub-domain
$D(s',r')\times \Cal O'$ of $D(s,r)\times \Cal O$ is the unique of
the homological equations (4.31-34), then $F$ satisfies the
skew-symmetric conditions:
$$F(-x,y,z,\bar z;\xi)=-F(x,y,z,\bar z;\xi),\ \forall \ (x,y,z,\bar z;\xi)
\in D(s',r')\times \Cal O'. \tag 4.40$$
\endproclaim
\demo{Proof} By (4.29,30) and (4.38), we see that $-F^z(-x)$
solves (4.38-2). Thus, $$-F^z(-x)=F^{\bar z}(x)$$ similarly, we
can show that $F^{x}(-x)=-F^x(x)$, $F^y(-x)=-F^y$. Consequently
$F(-x,y,z,\bar z;\xi)=-F(x,y,z,\bar z;\xi)$.$\qed$
\enddemo
%%%%%%%%%%%%%%
 \vskip 0.2cm

\noindent{\it 4.4. Solutions of the homological equations.}

\proclaim{Proposition 1} ( Solution of (4.36).) There is a subset
$\Cal O_+^1\subset\Cal O$ with $\text{Meas} \Cal
O_+^1\ge(\text{Meas} \Cal O)(1-K^{-1})$ such that for $\xi\in\Cal
O_+^1$,
$$|(k,\omega(\xi))|\gtrdot \frac{1}{K^{n+1}},
\text{ for all } 0\neq k\in\Bbb Z^n ,|k|\le K. \tag 4.41$$ Then,
on $D(s^\prime)\times \Cal O_+^1$, the equation (4.36) has a
solution $F^x(x,\xi)$ which is analytic in $x\in D(s^\prime)$ for
$\xi$ fixed and analytic in each $\xi_j,(j=1,..., n)$ for other
variables fixed, and which is real for real argument, such that
$$|\p_xF^x|_{D(s^\prime)\times\Cal O_+^1}\lessdot
K^Cr^2\varepsilon, \quad |\p_xF^x|_{D(s^\prime)\times\Cal
O_+^1}^{\Cal L}\lessdot K^Cr^2 \varepsilon^{\Cal L}.\tag 4.42$$
\endproclaim
\demo{Proof}  The existence of the set $\Cal O_+$ is well-known in
KAM theory. We omit the proof of the existence.  Recall $X_{\Cal
R}:D(s,r)\subset\Cal P\to \Cal P$ is real analytic in $(x,y,z,\bar
z)\in D(s,r)$ and each entry of $\xi\in\Cal O$. Expanding
$\p_x\Cal R^x$ into Fourier series
$$\p_x \Cal R^x=\sum_{0\neq k\in\Bbb Z^n,|k|\le K}\widehat{\p_x
\Cal R^x}(k)e^{\sqrt{-1}(k,x)}.\tag 4.43$$ Since $\p_x \Cal R^x$
is analytic in $x\in D(s)$, we get that the Fourier coefficients
$\widehat{\p_x \Cal R^x}(k)$'s decay exponentially in $k$, that
is,
$$|\widehat{\p_x \Cal R^x}(k)|\lessdot |\p_x
\Cal R^x|_{D(s)\times\Cal O}e^{-s|k|}\lessdot
e^{-s|k|}Kr^2\varepsilon ,\tag 4.44$$ where we have used (4.17).
Expanding $\p_x F^x$ into Fourier series:
$$F^x=\sum_{0\neq k\in\Bbb Z^n,|k|\le K}\widehat{\p_x
F^x}(k)e^{\sqrt{-1}(k,x)}.\tag 4.45$$ Inserting (4.43,45) into
(4.36), we get
$$\p_x F^x(x,\xi)=\sum_{0\neq k\in\Bbb Z^n,|k|\le K}\frac{\widehat{\p_x
\Cal R^x}(k)}{\sqrt{-1}(k,\omega)}e^{\sqrt{-}(k,x)}.$$ By
(4.36,41) as well as Lemma A.1, we get that for $x\in
D(s^\prime)$,
$$|\p_x F^x(x,\xi)|\lessdot K^Cr^2\varepsilon
\sum_{ k\in\Bbb Z^n}e^{-|k|(s-s^\prime)}\lessdot K^{\tilde
C}r^2\varepsilon .$$ Applying $\p_{\xi_j}$ to both sides of (4.36)
and using a method similar to the above, we can get the second
inequality of (4.42).$\qed$
\enddemo

\proclaim{Proposition 2} ( Solution of (4.37).) On
$D(s^\prime)\times \Cal O_+^1$,  the equation (4.37) has a
solution $F^y(x,\xi)$ which is analytic in $x\in D(s^\prime)$ for
$\xi$ fixed and analytic in each $\xi_j,(j=1,..., n)$ for other
variables fixed, and which is real for real argument, such that
$$|\p_xF^y|_{D(s^\prime)\times\Cal O_+^1}\lessdot
K^Cr\varepsilon, \quad |\p_xF^y|_{D(s^\prime)\times\Cal
O_+^1}^{\Cal L}\lessdot K^Cr\varepsilon^{\Cal L}.$$
\endproclaim
\demo{Proof} The proof follows almost exactly the proof of Prop.
1. We omit it.$\qed$
\enddemo
%%%%%%%%%%%%%%%%%%%%%%%%%%%
%%%%05/5/10下午3：45修改到此
%%%%%%%%%

 \subhead 5. Solution of (4.38)\endsubhead

This section is essential part of the present paper.   Let
$B_{ij}^{zz}$ ($i,j\in\zd$) be the elements of matrix $B^{zz}$.
Note that $B_{ij}^{zz}$ is a function of $(x,\xi)$. We temporally
omit the parameter $\xi$ for simplifying notation. Regarding
$B_{ij}^{zz}$ as a function of $x$, we let
$\widehat{B_{ij}^{zz}}(k)$ be the $k$-Fourier coefficient of
$B_{ij}^{zz}$. Let
$$\widehat{B_{ij}^{zz}}=(\widehat{B_{ij}^{zz}}(k-l):
|k|,|l|\le K).$$ Then $\widehat{B_{ij}^{zz}}$ is a matrix of order
$K^\sharp$. Let $\widehat{B}$ be a block matrix whose elements are
$\widehat{B_{ij}^{zz}}$ ($i,j\in \zd$), that is,
$$\widehat{B^{zz}}=(\widehat{B_{ij}^{zz}}:i,j\in \zd)=
(\widehat{B_{ij}^{zz}}(k-l):|k|,|l|\le K;i,j\in \zd).$$ Similarly
we have $\widehat{B^{z\bar z}}$ and $\widehat{B^{\bar z\bar z}}$.
Write $F^z$ as a column vector of infinite dimension:
$F^z=(F_j^z\in \Bbb C:j\in\zd)$. Let
$\widehat{F^z_j}=(\widehat{F_j^z}(k):|k|\le K)$ and
$\widehat{F^z}=(\widehat{F^z_j}:j\in\zd)$. Similarly we have
$\widehat{F^{\bar z}}$, $\widehat{R^{\bar z}}$ and
$\widehat{R^{\bar z}}$. Set
$$\widehat{F}=\pmatrix \widehat{F^z}\\ \widehat{F^{\bar
z}}\endpmatrix,\ \widehat{R}=\pmatrix -\sqrt{-1}\widehat{R^z}\\
\sqrt{-1}\widehat{R^{\bar z}}\endpmatrix\ \widehat{B}= \pmatrix
-\widehat{B^{z\bar z}}& \widehat{B^{z z}}\\ \widehat{B^{\bar z\bar
z}}&-\widehat{B^{z\bar z}}\endpmatrix.$$ Introduce a tensor
product space as follows
$$\un\ell^p:=\ell^p\otimes\Bbb C^{K^\sharp}=\{z=(z_j\in\Bbb C^
{K^{2n}}:j\in\zd): (||z||_p^*)^2=\sum_j|z|_j^2<\infty \}.$$ By
abuse use of notation, we $||\cdot||_p^*$ as $||\cdot||_p$.
 For any $u,v\in\un\ell^p$, we can regard
$(u, v)^T$ as a vector in $\un\ell^p\oplus\un\ell^p$.  Define
$$||(u,v)^T||_p=\sqrt{||u||_p+||v||_p}.
$$
Therefore,
$$||\widehat{F^u}||_p^2=\sum_{j\in\zd}\sum_{|k|\le
K}|\widehat{F^u_j}(k)|^2|j|^{2p}=\sum_{|k|\le
K}||\widehat{F^z}(k)||_p^2,\ u\in\{z,\bar z\},$$ where
$\widehat{F^u}(k) $ is $k$-Fourier coefficient of $F^u(x)$.
Similarly,
$$ ||\widehat{R^{
u}}||_p^2=\sum_{|k|\le K}||\widehat{R^{u}}(k)||_p^2,\ u\in\{z,\bar
z\}. $$ And
$$||\widehat
F||_p^2=||\widehat{F^z}||_p^2+||\widehat{F^{\bar z}}||_p^2.$$ In
view of (4.29,30), we see that the linear operator $\widehat{B}$
is self-adjoint in $\ell^0\oplus\ell^0$, i.e.,
$\widehat{B}=\widehat{B}^*$. With these notations, we have the
following lemma.

\proclaim{Lemma 5.1} (i). The vector  $\widehat{R}$ is in
$\un\ell^p\oplus\un\ell^p$ with
$$||\widehat{R}||_p\le Kr\varepsilon,\ ||\widehat{R}||_p^{\Cal L}
\le Kr\varepsilon^{\Cal L}\tag 5.1$$ (ii). The linear operator
$\widehat{B}$ is self-adjoint in $\un\ell^p\oplus\un\ell^p$ with
$p=0$, i.e., $\widehat{B}=\widehat{B}^*$; and
$$|||\widehat{B}|||_p, |||\widehat{B}|||_p^{\Cal L}\le\e\tag 5.2$$
\endproclaim
\demo{Proof} The conclusion (i) comes from (4.19). It follows form
(4.30) that the conclusion (ii) holds true.$\qed$
\enddemo

 Let
$$\Lambda_{\pm}=\text{diag}\ (\pm(k,\omega)+\Omega_j^0:|k|\le K,j\in
\zd),\ \Lambda=\pmatrix \Lambda_+& 0\\0&\Lambda_-\endpmatrix.
$$ Expand $\widehat{B^{zz}}, \widehat{B^{z\bar z}},\widehat{B^{\bar z\bar
z}}$, $F^z(x)$, $F^{\bar z}(x)$,  $R^z(x)$ and $R^{\bar z}(x)$ in
(4.38) into Fourier series,  we get a lattice equation which reads
$$(\Lambda+\widehat {B})\widehat{F}=\widehat{R}.\tag 5.3$$
Let $$M=(1+10K\max_{\xi\in\Cal O}{|\omega|})^{1/\kappa}.\tag 5.4$$
Write $\widehat{B}=(B_{ij}(k,l):|k|,|l|\le K, i,j\in\zd)$. Then we
see that $B_{ij}(k,l)$ is one of $\widehat{B^{zz}_{ij}}(k-l)$,
$\widehat{B^{z\bar z}_{ij}}(k-l)$ and $\widehat{B^{\bar z\bar
z}_{ij}}(k-l)$. Let
$$ \aligned &
\widehat{B}_{11}=(B_{ij}(k,l):|k|,|l|\le K, |i|,|j|\le M)\\
&\widehat{B}_{12}=(B_{ij}(k,l):|k|,|l|\le K, |i|\le M,|j|>M)\\&
\widehat{B}_{21}=(B_{ij}(k,l):|k|,|l|\le K, |i|>M,|j|\le M)\\&
\widehat{B}_{22}=(B_{ij}(k,l):|k|,|l|\le K, |i|>M,|j|>M
).\endaligned$$ Then
$$\widehat{B}=\pmatrix \widehat{B}_{11}&\widehat{B}_{12}
\\ \widehat{B}_{21} &\widehat{B}_{22}\endpmatrix.$$
And by (5.2),
$$|||\widehat{B}_{ij}|||_p, |||\widehat{B}_{ij}|||_p^{\Cal L}\le\e
,\ i,j\in\{1,2\}.\tag 5.5$$ Let
$$\aligned & \Lambda_{1}=\text{diag }(\pm(k,\omega)+\Omega^0_j:
|k|\le K, |j|\le M)\\&\Lambda_{2}=\text{diag
}(\pm(k,\omega)+\Omega^0_j: |k|\le K, |j|> M).\endaligned$$ Then
$$\Lambda=\pmatrix \Lambda_1&0\\0&\Lambda_2\endpmatrix.$$
In view of (5.4),
$$|\pm(k,\omega)+\Omega^0_j|\ge K, \forall\ |j|>M.\tag 5.6$$ Thus
by (5.5) and (5.6), we see that there exists the inverse of
$\Lambda_2+\widehat{B}_{22}$ and
$$|||(\Lambda_2+\widehat{B}_{22})^{-1}|||_p\le |||\sum_{j=0}^\infty
|||(\Lambda_2^{-1}\widehat{B}_{22})^j|||_p|||\Lambda_2^{-1}|||_p\lessdot
1.\tag 5.7$$ Set
$$\widehat B_{lef}:=\pmatrix E_1 & 0\\ -
(\varLambda_2+\widehat {B}_{22})^{-1}\widehat{B}_{21}& E_2
 \endpmatrix\tag 5.8$$
$$\Cal B_{rig}:=\pmatrix E_1 & -\widehat B_{12}(\varLambda_2+\widehat B_{22})^{-1}
\\
0&E_2 \endpmatrix,\tag 5.9$$ where $E_1$($E_2$, respectively) is a
unit matrix of the same order as that of $\Lambda_1$ ($\Lambda_2$,
respectively).  It is easy to verify that
$$|||\widehat B_{lef}^{-1}|||_p,\ |||\widehat B_{rig}^{-1}|||_p\lessdot 1.
\tag 5.10$$ Let
$$\widehat{ \Cal B}_{11}=
 \widehat{B}_{11}-\widehat B_{12}(\varLambda_2+\widehat B_{22})^{-1}
 \widehat B_{21}
 \tag 5.11$$
 Then there exists the inverse of $\varLambda+\widehat B $:
$$(\varLambda+\widehat B)^{-1} =\widehat B_{rig}^{-1}\pmatrix
(\varLambda_1+\widehat {\Cal B}_{11})^{-1}&0\\0& (\varLambda_2+
\widehat{
 B}_{22} )^{-1}\endpmatrix \widehat B_{lef}^{-1}\tag 5.12$$
provided that there exists the inverse of $\varLambda_1+\widehat
B_{11}$. And
$$|||(\varLambda+\widehat B)^{-1}|||_p^*\lessdot
|||(\varLambda_1+\widehat{\Cal B}_{11})^{-1}|||_p^*, \tag 5.13$$
where $*\equiv$ the blank or $\Cal L$. Now we are in position to
investigate the inverse of $\varLambda_1+\widehat{\Cal B}_{11}$.
First of all, we would like to point out that the matrix
$\widehat{\Cal B}_{11}$ is of finite order, and the order is
bounded by
$$K^*:=K^{2n} M^d\lessdot K^{2n+d\kappa^{-1}}.\tag 5.14 $$
Secondly, it follows from the self-adjoint-ness of $\widehat B$
that the matrix $\widehat{\Cal B}_{11}$ is also self-adjoint.
Thirdly, by (5.5,7,11) we have
$$|||\widehat{\Cal B}_{11}|||_p^{\Cal L}\lessdot \e .\tag 5.15 $$
Fourthly, since each element of $\widehat B$ is analytic in each
entry of $\xi\in\Cal O$, the matrix $\widehat{\Cal B}_{11}$ is
also analytic in each entry of $\xi\in\Cal O$. Without loss of
generality, we assume the first entry $\omega_1$ of $\omega$ is in
the interval $[1,2]$. Then the matrix $\Lambda_1+\widehat{\Cal
B}_{11}$ is non-singular if and only if
$\omega^{-1}_1\Lambda_1+\omega^{-1}_1\widehat{\Cal B}_{11}$ is
non-singular, and
$$|||(\Lambda_1+\widehat{\Cal
B}_{11})^{-1}|||_p\lessdot
|||(\omega^{-1}_1\Lambda_1+\omega^{-1}_1\widehat{\Cal
B}_{11})^{-1} |||_p. \tag 5.16
$$ Let
$$  \Xi \ \ :\sigma_1=1/\omega_1,\ \sigma_2=\omega_2/\omega_1,..., \sigma_n=
\omega_n/\omega_1. \tag 5.17$$ Then it is easy to get
$$\left|\text{det }\frac{\p (\sigma_1,...,\sigma_n)}{\p (\omega_1,...,
\omega_n)}  \right|=\omega_1^{-(n+1)}\gtrdot 1.\tag 5.18 $$
Combining (5.18) and (4.31), we can regard
$\sigma=(\sigma_1,...,\sigma_n)$ as a parameter vector instead of
$\xi$, and we can show that
$$\text{Meas }\Cal O\lessdot \text{Meas }\Xi(\Cal O)\lessdot
\text{Meas }\Cal O\tag 5.19 $$ as well as
$$|||\p_{\sigma}\widehat{\Cal
B}_{11}(\xi(\sigma))|||\lessdot \e,\tag 5.20 $$ $$
\p_\sigma\Omega_j^0=o(1),\tag 5.21$$ where (2.17) is used in
(5.21). Let
$$\Lambda_\sigma =\text{diag }(\pm(k_1+\sum_{l=2}^nk_l\sigma_l)+\sigma_1\Omega_j^0
: k=(k_1,...,k_n)\in\Bbb Z^n,|k|\le K,|j|\le M). \tag 5.22$$ Take
$\sigma$ as a parameter vector.  Then
$$\omega_1^{-1}\Lambda_1+\omega_1^{-1}\widehat{\Cal B}_{11}(\xi)=\Lambda_\sigma+
\sigma_1 \widehat{\Cal B}_{11}(\xi(\sigma)).\tag 5.23$$ Using
Lemma A.2, there are scalar functions
$\mu_1(\sigma),\cdots,\mu_{K^*}(\sigma)$ and a matrix-valued
function $U(\sigma)$ of order $K^*$ (See (4.14) for $K^*$), which
are analytic for each entry of real $\sigma$ and possess the
following properties for every $\sigma\in\Xi(\Cal O)$:
$$\varLambda_\sigma+\sigma_1\widehat{\Cal B}_{11}(\xi(\sigma))=U(\sigma)\text{diag}(\mu_1(\sigma),
\cdots, \mu_{K^*}(\sigma)) U^*(\sigma),\tag 5.24$$  and $$
U(\sigma)( U(\sigma))^*=E. \tag 5.25$$ It follows from (5.17) that
\footnote{ Here $|||\cdot|||$ is the $\ell_2$ norm of matrix.}
$$|||U(\xi)|||=|||U^*|||=1.\tag 5.26$$ Arbitrarily take
$\mu=\mu(\sigma)\in\{\mu_1(\sigma),...,\mu_{K^*}(\sigma)\}$. Let
$\phi$ be the normalized eigenvector corresponding $\mu$. Using
Lemma A.3 and noting (5.20) and (5.21), we get
$$\aligned \p_{\sigma_1}\mu=&((\p_{\sigma_1}(\varLambda_\sigma+
\sigma_1\widehat{\Cal B_{11}}(\xi(\sigma)))\phi, \phi)\\ =&
(\text{diag}(\Omega^0_j+\sigma_1\p_{\sigma_1}\Omega^0_j:|j|\le
M,|k|\le K )\phi,\phi)+\p_{\sigma_1} (\sigma_1\widehat
B(\xi(\sigma)))\phi,
\phi)\\=&\min_{j}\Omega_j^0+o(1).\endaligned\tag 5.27$$ Thus,
$$ \p_{\sigma_1}\mu\ge c>0.\tag 5.28$$
Let
$$ \Xi(\un{\Cal O}_l)=\{\sigma\in\Xi(\Cal O): |\mu_l|<
1/(KK^{*})\},\ l=1,...,K^*.\tag 5.29$$ By (5.28),
Meas$\Xi(\un{\Cal O}_l)<1/(KK^{*})$. Thus,
$$\text{Meas }\bigcup_{l=1}^{K^*}\Xi(\un{\Cal O}_l)<1/K.\tag 5.31$$
Let
$$\Xi(\un{\Cal O})=\Xi(\Cal O)\setminus \bigcup_{l=1}^{K^*}\Xi(\un
{\Cal O}_l).\tag 5.32$$ Therefore,
$$\text{Meas }\Xi(\un{\Cal O})\ge \text{Meas
}\Xi(\Cal O)(1-O(\frac{1}{K}))\tag 5.33$$ and for any
$\sigma\in\Xi(\un{\Cal O})$,
$$|\mu_l(\sigma)|\ge 1/KK^*.\tag 5.34$$ Moreover,
$$|||(\Lambda_\sigma+\sigma_1\widehat{\Cal B}_{11})^{-1}|||\le
|||U||||||U^*|||\max_l\{|\mu_l^{-1}|\}\le KK^*\le K^{C}\tag 5.35$$
where $C=n+d\kappa^{-1}+1$. Returning to the parameter $\xi$, by
(5.19),  there is a subset $\Cal O_+^2\subset\Cal O$ with
$$\text{Meas }\Cal O_+^2\ge (\text{Meas }\Cal O)(1-C_1K^{-1}),\tag
5.36$$ and for any $\xi\in\Cal O_+^2$
$$|||(\Lambda_1+\widehat{\Cal B}_{11})^{-1}|||\le C_2K^C.\tag 5.37$$
Observe that $\mu_l=\mu_l(\xi)$ is continuous in $\xi$. By means
of adjusting the constants $C_1, C_2$ in (5.36,37), we can extend
the set $\Cal O_+^2$ to an open set such that (5.36,37) still hold
true for any $\xi$ in the open set. Still denote by $\Cal O_+^2$
the open set.  Replacing the $\ell_2$ norm $|||\cdot|||$ by
$|||\cdot|||_p$ in (5.37), in view of the order of the matrix
$\Lambda_1+\widehat{\Cal B}_{11}$ being bounded by $K^*$, we get
$$|||(\Lambda_1+\widehat{\Cal B}_{11})^{-1}|||_p\le K^C(K^*)^p\le K^{\tilde C}
\tag 5.38$$ where $\tilde C$ is a constant depending
$n,d,\kappa,p$. By (5.13),
$$|||(\Lambda+\widehat{\Cal B})^{-1}|||_p\le  K^{\tilde C},\
\xi\in\Cal O_+^2. \tag 5.39$$ Note
$|||\p_\xi(\Lambda_1+\widehat{\Cal B})|||_p\le K$. We have
$$|||\p_\xi(\Lambda+\widehat{\Cal B})^{-1}|||_p\le  K^{1+\tilde C},\
\xi\in\Cal O_+^2. \tag 5.40$$%%%05年5月11日写到这里
By (5.1) and (5.39,40),
$$||\widehat F||_p\le K^{\bar C}r\varepsilon,\ \
||\p_\xi \widehat F||_p\le K^{\bar C}r\varepsilon^{\Cal L},\tag
5.41$$ where $\bar C$ is  a constant depending $p,d,n,\kappa$.
Recall
$$||\widehat F||_p^2=\sum_{|k|\le K}||\widehat{F^z}(k)||_p^2+
\sum_{|k|\le K}||\widehat{F^{\bar z}}(k)||_p^2.\tag 5.42$$ Recall
$$s^\prime=\frac{4}{K}.\tag 5.43$$
%%???? 5月14日：$s^\prime$ is too small so that $e^{|k|s^\prime}<\varepsilon$
%does NOT holds! ????? 5月15日：下一步的K应该取得更大才好。事实上，
% $K_m=2^mK_{m-1}$, $S_m=1/K_m$. %就很好！
Then for $u\in\{z,\bar z\}$,
$$\aligned \sup_{D(s^\prime)\times\Cal O_+^2}||F^u(x,\xi)||_p^2&=
\sup_{D(s^\prime)\times\Cal O_+^2}||\sum_{|k|\le
K}\widehat{F^u}(k)e^{\sqrt{-1}(k,x)}||_p^2\\& \le K\sum_{|k|\le
K}||\widehat{F^u}(k)||_p^2e^{2|k|s^\prime}\\&\le Ke^8\sum_{|k|\le
K}||\widehat{F^u}(k)||_p^2\\&\lessdot e^8K||\widehat F||_p^2\le
(K^Cr\varepsilon)^2
\endaligned\tag 5.44$$ That is,
$$\sup_{D(s^\prime)\times\Cal O_+^2}||F^u(x,\xi)||_p\lessdot K^Cr\varepsilon
,\ u\in\{z,\bar z\}.\tag 5.45$$ Similarly,
$$\sup_{D(s^\prime)\times\Cal O_+^2}||\p_\xi F^u(x,\xi)||_p\lessdot K^Cr\varepsilon
,\ u\in\{z,\bar z\}.\tag 5.46$$ Consequently, we have the
following proposition. \proclaim{Proposition 3} There is a subset
$\Cal O_+^3\subset\Cal O$ with $\text{Meas }\Cal O_+^3=
(\text{Meas }\Cal O)(1-O(1/K))$ and there is  functions $F^z,
F^{\bar z}: D(s^\prime)\times\Cal O_+^3\to\ell^p$ which are
analytic in $x\in D(s^\prime)$ and also analytic in each entry of
$\xi\in\Cal O_+^3$; moreover the functions $F^z, F^{\bar z}$ solve
(4.38) and satisfy the estimates (5.45,46).
\endproclaim

\subhead 6.  Symplectic change of variables
\endsubhead

In this section, our procedure is standard and almost the same as
that of Section 3 in [P1,p.128-132]. Here we give out the outline
of the procedure. See [P1] for the details. In this section, we
denote by $C$ a universal constant depending only on
$n,d,p,\kappa$. The constant $C$ might be different in different
places.

 {\it 6.1. Coordinate transformation.}
By Propositions 1-3, we get a Hamiltonian  $F$ on $D(s^\prime,r)$
where
$$F=F^x+(F^y,y)+\<F^z,z\>+\<F^{\bar z},z\>$$ and give estimates of
$F^x$, $F^y$, $F^z$ and $F^{\bar z}$. Let $X_F$ be the vector
field corresponding to the Hamiltonian $F$, that is,
$$X_F=(-\p_y F,\p_x F,\sqrt{-1}\p_z F, -\sqrt{-1}\p_{\bar z}F),$$
here $\p_z $  and $\p{\bar z}$ are the usual $\ell^2$-gradients.
Recall $s=4/K_-$ and $s^\prime=4/K$. And let
$s^{\prime\prime}=3/K$, $s^{\prime\prime\prime}=2/K,
s^{\prime\prime\prime\prime}=1/K$. Let
$$\Cal O_+=\bigcup_{j=1}^2\Cal O_+^j.$$
 It follows from Prop.1,2 and 3 that for  $(x,y,z,\bar z;\xi)\in
 D(s^\prime,r)\times\in\Cal O_+$,
$$\aligned \ _r\W X_F\W_{p}=&|\p_y F|+\frac{1}{r^2}|\p_x F|+
\frac{1}{r}||\p_zF|| _{ p}+\frac{1}{r}||\p_{\bar z}F|| _{ p}\\=&
|F^y|+\frac{1}{r^2}|\p_x F^x|+\frac{1}{r} ||F^z||_{ p}+\frac{1}{r}
||F^{\bar z}||_{ p}\\ \le & |F^y|+\frac{K}{r^2}| F^x|+\frac{1}{r}
||F^z||_{ p}+\frac{1}{r} ||F^{\bar z}||_{ p}\\ \lessdot & K^C\cdot
\varepsilon,
\endaligned $$ where $C$ is a large constant depending only on $n,\kappa,
d$ and $p$. That is,
$$\ _r\W X_F\W_{p, D(s^\prime,r)\times\Cal O_+}\lessdot K^C
 \varepsilon.\tag 6.1$$
 Similarly, we have
$$\ _r\W X_F\W_{p, D(s^{\prime},r)\times\Cal O_+}^{\Cal L}\lessdot
K^C\varepsilon^{\Cal L}. \tag 6.2$$ where we have used the
assumption $\varepsilon<\varepsilon^{\Cal L}$.  As in [P1,p.129],
we introduce the operator norm
$$_r\W\W L\W\W_{p}=\sup_{W\neq 0}\frac{_{r}\W LW\W_{p}}
{_r\W W\W_{p}}. \tag 6.3$$  Using (6.1), (6.2) and the generalized
Cauchy's inequality (See Lemma A.3 of [P1,p.147]) and the
observation that every point in $D(s^{\prime\prime},r/2)$  has at
least $_r\W\cdot\W_{p}$ distance $1/9K$ to the boundary of
$D(s^\prime,r)$, we get
$$ \sup_{D(s^{\prime\prime},r/2)\times\Cal O_+}\ _r
\W\W DX_F\W\W_{p} \lessdot K\ _r\W X_F\W_{p,
D(s^\prime,r)\times\Cal O_+} \le K^C\varepsilon.\tag 6.4$$
$$ \sup_{D(s^{\prime\prime},r/2)\times\Cal O_+}\ _r\W\W
DX_F\W\W_{p}^{\Cal L} \lessdot K\ _r\W X_F\W_{p,
D(s^\prime,r)\times\Cal O_+}^{\Cal L} \le K^C\varepsilon^{\Cal
L},\tag 6.5$$ where $ DX_F$ is the differential of $X_F$.
 Assume that $K^C\varepsilon$ and
 $K^C\varepsilon^\Cal L$ are small
enough. (These assumptions will be fulfilled in the following KAM
iterations. In fact, in $m$-th KAM step, $K\approx C2^{m^2/2}$ and
$\varepsilon=\e^{(4/3)^m}$. It follows that $K^C\varepsilon\ll
1$.) Arbitrarily fix $\xi\in\Cal O_+$. By (6.1), the flow $X_F^t$
of the vector field $X_F$ exists on
$D(s^{\prime\prime\prime},r/4)$ for $t\in[-1,1]$ and takes the
domain into $D(s^{\prime\prime},r/2)$, and by Lemma A.4 of [P1,
p.147], we have
$$_r\W X_F^t-id\W_{p,D(s^{\prime\prime\prime},r/4)\times\Cal O_+}\lessdot
\ _r\W X_F\W_{p,D(s^\prime,r)\times\Cal O_+}\lessdot
K^C\varepsilon\tag 6.6
$$ and
$$\aligned & _r\W X_F^t-id\W_{p,D(s^{\prime\prime\prime},r/4)\times
\Cal O_+} ^{\Cal L}\\ &\lessdot \exp( _r\W\W
DX_F\W\W_{p,D(s^{\prime\prime},r/2)\times\Cal O_+})\cdot\  _r\W
X_F\W_{p,D(s^\prime,r)\times\Cal O_+}^{\Cal L}\\&\lessdot
(\exp(K^C\varepsilon))K^C\varepsilon^{\Cal L}\lessdot
K^C\varepsilon^{\Cal L},
\endaligned\tag 6.7
$$ for $t\in [-1,1]$. Furthermore, by the generalized Cauchy's
inequality, $$_r\W\W
DX_F^t-I\W\W_{p,D(s^{\prime\prime\prime\prime},r/8)\times\Cal
O_+}\lessdot K^C\varepsilon,\ \ t\in[0,1]\tag 6.8$$ and
$$_r\W\W
DX_F^t-I\W\W_{p,D(s^{\prime\prime\prime\prime},r/8)\times\Cal
O_+}^{\Cal L}\lessdot K^C\varepsilon^{\Cal L},\ \ t\in[0,1]\tag
6.9$$

{\it 6.2. The new error term.} Subjecting $H=N+\grave R$ to the
symplectic transformation $\Phi=X_F^t|_{t=1}$ we get the new
Hamiltonian scale $H_+:=H\circ \Phi=H\circ X_F^1$ on
$D(s^{\prime\prime\prime\prime},\eta r)$ where $0<\eta\ll 1$.
 By
Taylor's formula
$$\aligned H_+&=(N+R)\circ X_F^1=(N+R+(\grave R-R))\circ X_F^1\\&=
 (N+\Cal R+R_K+(\grave R-R))\circ X_F^1\\
&=N-\{F,N\}+\int_0^1\{t\{F,N\},F\}\circ X_F^t\ dt\\
&\quad +\Cal R+\int_0^1\{\Cal R,F\}\circ X_F^t\ dt+(R_K+(\grave
R-R))\circ X_F^1.
\endaligned\tag 6.10$$ Recall that (4.39) holds true when $F$ solves (4.36-38).
Thus,
$$H_+=N_++\grave R_+\tag 6.11$$
 where
$$N_+=N+\frac12\<(F^y,\p_xB^{zz})z, z\>+\<(F^y,
\p_xB^{z\bar z}) z,\bar z\> +\frac12\<(F^y,\p_xB^{\bar z\bar
z})\bar z,\bar z\>\tag 6.12$$
$$\grave R_+=\grave R_+^1+\grave R_+^2+ \grave R_+^3\tag 6.13$$
where $$\aligned \grave R_+^1=&
\sqrt{-1}\left\<(1-\Gamma_K)(B^{\bar z\bar z}F^z),\bar z\right\>+
\sqrt{-1}\left\<(1-\Gamma_K)(B^{z\bar z}F^{ z})
z\right\>\\&-\sqrt{-1}\left\<(1-\Gamma_K)(B^{ z z}F^{\bar z}),
z\right\>- \sqrt{-1}\left\<(1-\Gamma_K)(B^{z\bar z}F^{\bar z}),
\bar z\right\>\endaligned\tag 6.14$$ $$\grave R_+^2=R_K\circ
X_F^1+(\grave R-R)\circ X_F^1\tag 5.15$$ $$\grave
R_+^3=\int_0^t\{(1+t)\Cal R,F\} \circ X_F^t\ dt.\tag 6.16$$ By
(4.26,27,28),
$$N_+=(\omega_+,y)+\<\Omega^0z,\bar z\>+\frac12\<B^{zz}_+z,z\>
+\<B^{z\bar z}_+z,\bar z\>+\frac12\<B^{\bar z\bar z}_+\bar z,\bar
z\>\tag 6.17$$ with
$$\aligned & B^{zz}_+=\un B^{zz}+R^{zz}+(F^y,\p_x(\un B^{zz}+R^{zz})),
\\& B^{z\bar z}_+=\un B^{z\bar z}+R^{zz}+(F^y,\p_x(\un B^{z\bar z}+
R^{z\bar z})),\\
& B^{\bar z\bar z}_+=\un B^{\bar z\bar z}+R^{\bar z\bar
z}+(F^y,\p_x(\un B^{zz}+R^{\bar z\bar z})),\endaligned\tag 6.18$$
and
$$\omega_+=\omega=\un\omega+\widehat{F^y}(0).\tag 6.19$$
Hence, the new perturbing vector field is
$$X_{\grave R_+}=X_{\grave R^1_+}+(X_F^1)^*(X_{\grave R}-X_R+R_K)+\int_0^t(X_F^t)^*
[X_{(1+t)\Cal R}, X_F]\ dt,$$ where $(X_F^t)^*$ is the pull-back
of $X_F^t$ and $[\cdot,\cdot]$ is the commutator of vector fields.
 We are now in position to estimate the new perturbing vector field
 $X_{\grave R_+}$. Let $Y:D(s^\prime,r)\subset\Cal P\to \Cal P$ be a vector field on
 $D(s^\prime,r)$, depending on the parameter
 $\xi\in \Cal O_+$. Let $U=D(s^{\prime\prime\prime\prime},\eta r)
 \times\Cal O_+$ and
 $V=D(s^{\prime\prime\prime},2\eta r)\times\Cal O_+$ and
 $W=D(s^{\prime\prime},4\eta r)\times\Cal O_+$. By (6.6) and the ``proof of estimate $(12)$" of [P1,
 p.131-132]\footnote{Let $a=0$ in (12) of [P1].}, we have that
$$_{\eta r}\W (X_F^t)^* Y\W_{p, U}\lessdot\  _{\eta r}\W Y\W_{p, V}
\tag 6.20$$ and
$$_{\eta r}\W (X_F^t)^* Y\W_{p, U}^{\Cal L}\lessdot\  _{\eta r}\W
Y\W_{p, V}^{\Cal L}+\frac{K}{\eta^2}\ _{\eta r}\W Y\W_{p, W}\cdot\
_{\eta r}\W X_F\W_{p, V}^{\Cal L}.\tag 6.21$$ We assume that
$$\varepsilon K^C/\eta^2\lessdot 1.\tag 6.22$$ These assumptions will be
fulfilled in the KAM iterative lemma later. (In fact, in $m$-th
KAM iterative step, we will let $K^C\approx 2^{Cm^2},
\varepsilon=\e^{2^m}, \eta=\varepsilon^{1/3}$. This implies (6.22)
is fulfilled.) By (4.10) and (6.20,21),
$$ _{\eta r}\W (X_F^1)^*(X_{\grave R}-X_R)\W_{p, U}\lessdot
\ _{\eta r}\W X_{\grave R}-X_R\W_{p, V}\lessdot
\eta\varepsilon\tag 6.23$$ and
$$_{\eta r}\W (X_F^1)^*(X_{\grave R}-X_R)\W_{p, U}^{\Cal L}
\lessdot \eta\varepsilon^{\Cal L}+\frac{\varepsilon K^C \eta
}{\eta^2} \varepsilon^{\Cal L}\lessdot\eta\varepsilon^{\Cal
L}.\tag 6.24$$
   Recall
  that (4.9) holds still true after replacing $R$ by $\Cal R$.
By (4.2) and (5.4,5) and using the generalized Cauchy estimate,
following [P1,p.130-131] we get
$$\aligned  _r\W [X_{\Cal R(t)},X_F]\W_{p,U }& \lessdot K
\ _r\W X_{\Cal R}\W_{p, V} \cdot \ _r\W X_F\W_{p,V}\\&\lessdot K \
_r\W X_{\grave R}\W_{p, W} \cdot \ _r\W X_F\W_{p,W} \\&\lessdot
K^{1+C}\varepsilon ^2<\eta\varepsilon\endaligned\tag 6.25$$ and
$$\aligned & _r\W [X_{\Cal R(t)},X_F]\W_{p,U }^\Cal L\\& \lessdot K
\ _r\W X_{\grave R}\W_{p, D_1}^\Cal L\ _r\W X_F\W_{p,W}+K\ _r\W
X_{\grave R}\W_{ p,
 W}\ _r\W
X_F\W_{p,W}^{\Cal L}\\& \lessdot K^C\varepsilon ^\Cal L
\varepsilon+K^C\varepsilon\varepsilon^\Cal L <\eta\varepsilon^\Cal
L\endaligned\tag 6.26$$ Finally, we have
$$_{\eta r}\W Y\W_{p,U }\lessdot  \eta^{-2}\ _r\W Y\W_{p,U
},\quad _{\eta r}\W Y\W_{p,U }^\Cal L\lessdot \eta^{-2}\ _r\W
Y\W_{p,U }^\Cal L,\tag 6.27 $$ for any vector field $Y$. Note that
for any function $f:D(s)\times\Cal O\to \ell^p$ or $\Bbb R^n$
which is analytic in $x\in D(s)$ and in each entry of $\xi\in\Cal
O$,
$$\aligned &||(1-\Gamma_K)f||_{D(s^\prime)\times\Cal O}\le
\sum_{|k|>K}||\widehat{f}(k)||_pe^{s^\prime|k|}\\& \le
||f||_{D(s)\times\Cal O}\sum_{|k|>K}e^{-(s-s^\prime)|k|}\\& \le
||f||_{D(s)\times\Cal O}\cdot\varepsilon^2,\endaligned\tag 6.28$$
where we have used
$$(s-s^\prime)K\approx 2|\ln\varepsilon|.\tag 6.29$$ Applying (6.27,28) to (6.14), we can easily get
$$ _{\eta r}\W X_{\grave R_+^1}\W_{p,U}\lessdot \eta^{-2}\varepsilon^2
\lessdot \eta\varepsilon,\ _{\eta r}\W X_{\grave
R_+^1}\W_{p,U}^{\Cal L}\lessdot \eta\varepsilon^{\Cal L}.\tag
6.30$$
 Collecting all terms above, we then arrive at the estimates
$$_{r_+}\W X_{\grave R_+}\W_{p,D(s_+,r_+)\times\Cal O_+ }\lessdot
\eta\varepsilon,\quad _{r_+}\W X_{\grave
R_+}\W_{p,D(s_+,r_+)\times \Cal O_+ }^\Cal L\lessdot
\eta\varepsilon^\Cal L. \tag 6.31$$ where
$$s_+=s^{\prime\prime\prime\prime},\quad r_+=\eta
r.\tag 6.32$$ It follows from (4.22) and (4.40) that the
perturbing Hamiltonian $\grave R_+$ satisfies the same symmetric
condition as $\grave R$, that is:
$$\grave R_+(-x,y,z,\bar z;\xi)= \grave R_+(x,y,\bar z,z;\xi),\ \ \forall\
(x,y,z,\bar z;\xi)\in D(s_+,r_+)\times \Cal O_+.\tag 6.33$$
Finally we can easily check that $N_+$ satisfies the same
conditions as $\un B$, that is, the conditions (4.2-6) hold true
by replacing $\un B$ by $B_+$ and replacing $\un\omega$ by
$\omega_+$. Here we omit the details.

%%%%%%%%%%%%%%
%%%%%%%%%%%%%%
%%%%%%%%%%%%%%
%%%%%%%%%%%%%%
%%%%%%%%%%%%%%
%%%%%%%%%%%%%%
%%%%%%%%%%%%%%
%%%%%%%%%%%%%%
%%%%%%%%%%%%%%
\subhead 7. Iterative lemma and proof of the theorem\endsubhead

\noindent{\it 7.1.  Iterative constants.} As usual, the KAM
theorem is proven by the Newton-type iteration procedure which
involves an infinite sequence of coordinate changes. In order to
make our iteration procedure run, we need the following iterative
constants:

\item{1.} $\e_0=\e$, $\e_l=\e_0^{\wedge}(4/3)^l$, $l=1,2,...$;

\item{2.} $\eta_l=\e_l^{1/3}$, $l=0,1,2,..$;

%\item{4.} $e_0=0$, $ e_m= (1^{-2}+\cdots+m^{-2})/2\sum_{j=1}^\infty j^{-2}$,
 %(thus, $e_m<1/2$ for all $m\in \Bbb N$);

\item{3.} $K_l=2^{\sum_{j=1}^l j}|\ln\e|$, $l=1,2,...$, (thus,
$K_l=2^{l(l+1)/2} |\ln\e|$); \item{4.} $s_0>0$ is given,
$s_l=1/K_{l-1}$, $l=1,2,...$, (thus, $(s_l-s_{l+1})K_l\ge
2|\ln\e_l|$.); \item{5.} $r_0>0$ is given, $r_l=\eta_l r_0$,
$l=1,2,...$; \item{6.} $D(s_l)=\{x\in\Bbb C^n/(2\pi)^n:|\Im
x|<s_l\}$ \item{7.} $D(s_l,r_l)=\{(x,y,z,\bar z)\in\Cal P:|\Im
x|<s_l,|y|<r_l^2,||z||_l< r_l, ||\bar z||_p<r_l\}$.
%\item{?.} $\Cal O_m$'s are the $\nu_m$-neighborhood of $\Cal O_m$ in $\Bbb R^n$.

\vskip0.3cm
%%06年5月16日写到这里
\noindent {\it 7.2.  Iterative Lemma.  Consider a family of
Hamiltonian functions $ H_l$ ($0\le l\le m$):
$$\aligned H_l=&(\omega_l(\xi),y)+(\Omega^0(\xi)z,\bar z)+\frac12\<
B^{zz}_l(x,\xi)z, z\>+\<B^{z\bar z}_l(x,\xi)z, \bar
z\>\\&+\frac12\<B_l^{\bar z\bar z}(x,\xi)\bar z,\bar z\>+\grave
R_l(x,y,z,\bar z,\xi)\endaligned \tag 7.1$$ where
$B_l^{zz}(x;\xi)$'s are analytic in $x\in D(s_l)$ for any fixed
$\xi\in\Cal O_l$, and all of $\omega_l(\xi),\Omega^0(\xi)$ and $
B_l^{zz}(x;\xi)$'s (for fixed $x\in D(s_l)$) are analytic in each
entry $\xi_l$ ($l=1,..., n$) of $\xi\in\Cal O_l$, and $\grave
R_l(x,y,z,\bar z,\xi)$ is analytic in $(x,y,z,\bar z)\in
D(s_l,r_l)$ and analytic in each entry of $\xi\in\Cal O_l$. Assume
$B_l^{zz}=(B^{zz}_{l,ij}:i,j\in\zd)$'s satisfy the following
conditions:

\par\noindent {\it ($l.1$). Symmetry.} For any $(x,\xi)\in D(s_l)\times \Cal O_l$,
the operators $B^{zz}_l, B^{z\bar z}$ and $ B_{\bar z\bar z}$ are
self-adjoint from $\ell_2$ to $\ell_2$.

\par\noindent {\it ($l.2$). Finiteness of Fourier modes.}
$$ B_l^{zz}(x;\xi)=\sum_{|k|\le K_l}\widehat{B_{l}^{zz}}(k)
e^{\sqrt{-1}(k,x)},\cdots $$
\par\noindent {\it ($l.3$). Boundedness.}
$$\sup_{D(s_l)\times\Cal O_l}|||B_l^{zz}(x;\xi)|||^*_p\lessdot  \e,\ ...$$
where $*\equiv$ the blank and $\Cal L$.
 Moreover, we assume the parameter sets $\Cal O_l$'s satisfy
\noindent $(l.4)$. $$\Cal O_0\supset\cdots\supset\cdots\Cal O
_l\supset\cdots\Cal O_m$$ with
$$\text{ Meas }\Cal O_l\ge (\text{ Meas
}\Cal O_0)(1-K_l^{-1});$$

 \noindent $(l.5)$.  The map $\xi\mapsto \omega_l(\xi)$ is
analytic in each entry of  $\xi\in\Cal O_l$,  and
$$\inf_{\Cal O_l}\left|\text{det }\frac{\p \omega_l}{\p \xi}\right|
\ge c_1>0,\ \sup_{\Cal O_l}|\p_\xi^j \omega_l|\le c_2,\ j=0,1.$$

\noindent $(l.6)$. The perturbation $\grave R_l(x,y,z,\bar z
;\xi)$ is analytic in the space coordinate domain $D(s_l,r_l)$ and
also analytic in each entry $\xi_k$ ($k=1,...,n$) of the parameter
vector $\xi\in\Cal O_l$, and is real for real argument; moreover,
its Hamiltonian vector field $X_{\grave R_l}:=(\p_y\grave
R_l,-\p_x\grave R_l,\sqrt{-1}\p_z \grave R_l,-\sqrt{-1}\p_{\bar z}
\grave R_l )^T$ defines on $D(s_l,r_l)$ a analytic map
$$X_{\grave R_l}:D(s_l,r_l)\subset\Cal P\to\Cal P. $$

\noindent $(l.7)$. In addition, the vector field $X_{\grave R_l}$
is analytic in the domain $D(s_l,r_l)$ with small norms
$$_{r_r}\W X_{\grave
R_l}\W_{p,D(s_l,r_l)\times\Cal O_l}\lessdot \e_l,\quad _{r_l}\W
X_{\grave R_l}\W_{p,D(s_l,r_l)\times\Cal O_l}^{\Cal L}\lessdot
\e_l^{1/3}.$$ Then there is is an absolute positive constant
$\epsilon^*$ enough small such that, if $0<\epsilon_0<\epsilon^*$,
there is a set $\Cal O_{m+1}\subset \Cal O_m$, and a change of
variables $\Phi_{m+1}: \Cal D_{m+1}:=D(s_{m+1},r_{m+1})\times \Cal
O_{m+1}\to D(s_{m},r_{m})$ being real\footnote{The word ``real"
means $\overline{\Phi_{m+1}(z,\xi)}=\Phi_{m+1}(\bar z,\xi)$ for
any $(z,\xi)\in\Cal D_{m+1}$.}, analytic in $$(x,y,z,\bar z)\in
D(s_{m+1},r_{m+1})$$ and each entry $\xi\in\Cal O_{m+1}$, as well
as following estimates holds true:
$$_{r_m}\W \Phi_{m+1}-id\W_{p,\Cal D_{m+1}}\lessdot
\e_{m}^{1/2}
$$ and
$$_{r_m}\W \Phi_{m+1}-id\W_{p,\Cal D_{m+1}}^{\Cal L}\lessdot
\e_{m}^{1/4}.$$  Furthermore, the new Hamiltonian
$H_{m+1}:=H_m\circ \Phi_{m+1}$ of the form
$$\aligned H_{m+1}=&(\omega_{m+1}(\xi),y)+(\Omega^0(\xi)z,\bar z)+
\frac12\< B^{zz}_{m+1}(x,\xi)z, z\>+\<B^{z\bar z}_{m+1}(x,\xi)z,
\bar z\>\\&+\frac12\<B_{m+1}^{\bar z\bar z}(x,\xi)\bar z,\bar
z\>+\grave R_{m+1}(x,y,z,\bar z,\xi)\endaligned \tag 7.2$$
  satisfies all the
above conditions $(l.1-7)$ with $l$ being replaced by $m+1$.}
\vskip 0.3cm
%%%%%%%%
%%%%%%%%

\noindent{\it 7.3. Proof of The Iterative Lemma.}

As stated as in the iterative lemma, we have got a family of
Hamiltonian functions $H_l$'s ($l=0,1,...,m$) which satisfy the
conditions $(l.1-7)$. We now consider the Hamiltonian $H_m$. Let
$\un N=N_m$, $\grave R=\Grave R_m$, $N_+=N_{m+1}$ and $\grave
R_+=\grave R_{m+1}$. Let $s=s_m$,
$s^{\prime\prime\prime\prime}=s_{m+1}$, $\eta=\eta_m$,
$r=r_m=\eta_m r_0$, $\varepsilon=\e_m$, $\varepsilon^\Cal
L=\e_m^{1/3}$ and $\Cal O_{m+1}=\Cal O_+^1\cup\Cal O_+^2$.
Clearly, $\varepsilon<\varepsilon^\Cal L$. By means of the
conclusion in Section 4, we got that there is a Hamiltonian
$F=F_m$ defined on\footnote{Note $\Cal
D_{m+1}:=D(s_{m+1},r_{m+1})\times\Cal O_{m+1}\subset
D(s_{m+1},r_{m+1})\times\Cal O_{m+1}$. }
$D(s_{m+1},r_{m+1})\times\Cal O_{m+1} $ and a symplectic change of
variables $\Phi_{m+1}=X_{F_m}^t|_{t=1}$. Note
$\eta_n\e_m=\e_{m+1}$. It is easy to verify that the conditions
($m+1.1-7$) are fulfilled. We omit the details.

\vskip 0.3cm

%%%%%%%%%%%%
%%%%%%%%%%%
%%%%%%%%%%%%

\noindent{\it 7. Proof of the theorem 2.1}

 The proof is similar
to that of [P1]. Here we give an outline.  By Assumptions A,B, C,
D, E  and the smallness assumption in Theorem 2.1, the conditions
$(l.1-7)$ in the iterative lemma in Section 6.2 are fulfilled with
$l=0$. Hence the iterative lemma applies to $\tilde H$.
Inductively, we get what as follows:

\item{(i).} Domains: for $m=0,1,2,...$,
$$ \Cal D_m:=D(s_{m},r_{m})\times \Cal O_{m},\quad \Cal D_{m+1}\subset\Cal D_m
 ;$$

\item{(ii).} Coordinate changes:
$$\Psi^m=\Phi_1\circ\cdots\circ\Phi_{m+1}:\Cal D_{m+1}
\to  D(s_0,r_0),;$$

\item{(iii).} Hamiltonian functions $\tilde H_m$  ($m=0,1,...$)
satisfy the conditions $(l.1,2,3)$ with $l$ replaced by $m$;

Let $\Cal O_\infty=\cap_{m=0}^\infty\Cal O_m$, $\Cal D_\infty=\cap
\Cal D_m$. By the same argument as  in [P1, pp.134], we conclude
that $\Psi^m,D\Psi^m,\tilde H_m,X_{H_m}$ converges uniformly on
the domain $\Cal D_\infty$, and $X_{\tilde
H_\infty}\circ\Psi^\infty=D\Psi^\infty\cdot X_{\omega_*}$ where
$$\aligned \tilde H_{\infty}:=&\lim_{m\to\infty}\tilde H_m=
(\omega_*(\xi),y)+\<\Omega^0(\xi)z,\bar z\>+\\&+ \frac12
\<B^{zz}_\infty(x,\xi) z,z\>+\frac12 \<B^{z\bar z}_\infty(x,\xi)
z,\bar z\>+\frac12 \<B^{\bar z\bar z}_\infty(x,\xi) \bar z,\bar
z\>
\endaligned$$ here
$B^{zz}_{\infty}=\lim_{m\to\infty}B_{m}^{zz}$,..., and
$X_{\omega_*}$ is the constant vector field $\omega_*$ on the
torus $\Bbb T^n$.
 Thus, $\Bbb T^N\times \{0\}\times\{0\}$ is an
embedding torus with rotational frequencies
 $\omega_*(\xi)\in \omega_*(\Cal O_\infty) $ of the Hamiltonian $\tilde H_\infty$.
 Returning the original Hamiltonian $\tilde H$, it has an
 embedding torus $\Phi^\infty (\Bbb T^n\times \{0\}\times\{0\} )$
 with  frequencies  $\omega_*(\xi)$. This proves the Theorem.$\qed$

\subhead 9. Appendix A. Some Technical lemmas\endsubhead

\proclaim{Lemma A.1} For $\mu>0,\nu>0$,the following inequality
holds true:
$$\sum_{k\in \Bbb Z^d}e^{-2|k|\mu}|k|^\nu\le (\frac{\nu}{e})^\nu\frac{1}{\mu^{\nu+d}}(1+e)^d.$$
\endproclaim
\demo{Proof} This Lemma can be found in [B-M-S].
\enddemo

\proclaim{Lemma A.2} Consider an $n\times n$ complex matrix
function $Y(\xi)$ which depends on the real parameter $\xi\in\Bbb
R$. Let $Y(\xi)$ be a matrix function satisfying conditions:

\item{(i)} $Y(\xi)$ is self-adjoint for every $\xi\in\Bbb R$;
i.e., $Y(\xi)=(Y(\xi))^*$, where star denotes the conjugate
transpose matrix;

\item{(ii)} $Y(\xi)$ is an analytic function of the real variable
$\xi$.

\noindent Then there exist scalar functions
$\mu_1(\xi),\cdots,\mu_n(\xi)$ and a matrix-valued function
$U(\xi)$, which are analytic for real $\xi$ and possess the
following properties for every $\xi\in\Bbb R$:
$$Y(\xi)=U(\xi)\text{diag}(\mu_1(\xi),\cdots,\mu_n(\xi))U^*(\xi),\quad U(\xi)( U(\xi))^*=E.$$
\endproclaim
\demo{Proof} See [pp.394-396, G-L-R]. $\qed$

It is worth to point out that this lemma does not hold true for
$\xi\in\Bbb R^k $ with $k>1$. See [Ka].
\enddemo
\proclaim{Lemma A.3} Assume $Y=Y(\xi)$ satisfies the conditions in
Lemma A.3. Let $\mu=\mu(\xi)$ be any eigenvalue of $Y$ and $\phi$
be the normalized eigenfunction corresponding to $\mu$. Then
$$\p_\xi\mu=((\p_\xi Y)\phi,\phi)).$$

\endproclaim
\demo{Proof} The proof can be found in [Ka,p.125].
\enddemo

  \proclaim {
References}\endproclaim

\item{[Ar]} Arnol$'$d, I. V.: Mathematical  Methods of Classical
Mechanics.  Springer-Verlag,  New York (1978)

\item{[B]} Bambusi, D., Birkhoff normal form for some nonlinear
PDEs, Commun. Math. Phys., {\bf 234}, 253--285 (2003)

\item{[B-B]} Berti, M. and Bambusi, D., A Birkhoff-Lewis type
theorem for some Hamiltonian PDEs, preprint, (2003) \item{[B-G]}
Bambusi, D. and Gr\'ebert, Birkhoff normal form for PDEs with TAME
modulus, Duke Math. {\it to appear}

\item{[B-P]} Bambusi, D., and Paleari, S., Families of periodic
orbits for resonant PDE's, J. Nonlinear Science, {\bf 11} (2001),
69-87
 \item{[Bo-K]} Bobenko, A. I., and  Kuksin, S.,
 The nonlinear Klein-Gordon equation on an interval as perturbed
    sine-Gordon equation,  Comm. Math. Helv.
    {\bf 70}, 63--112 (1995)

\item{[Br-K-S]} Bricmont, J., Kupiainen, A. and Schenkel, A.,
Renormalization group and the Melnikov Problems for PDE's, Commun.
Math. Phys., {\bf 221}, 101--140 (2001)

\item{[B-M-S]} Bogolyubov, N. N., Mitropolskij, Yu. A. and
Samojlenko, A. M.: Methods of Accelerated Convergence in Nonlinear
Mechanics. Springer-Verlag, New York (1976) [Russian Original:
Naukova Dumka, Kiev, 1969]

%\item{[Bi]} Birkhoff, G. D., Dynamical systems, Am. Math. Soc.
% Colloq. Publ. IX. New York: American Mathematical Society. VIII, 1927

\item{[Bo1]}   Bourgain, J., Construction of quasi-periodic
solutions for Hamiltonian
    perturbations of linear
     equations and application to nonlinear pde,
     Int. Math. Research Notices {\bf 11}  475--497(1994)
\item{[Bo2]}  Bourgain, J., Periodic solutions of nonlinear wave
equations, Harmonic analysis and partial equations, Chicago Univ.
Press, 1999, pp.69-97

 \item{[Bo3]}  Bourgain, J., Quasi-periodic solutions of
Hamiltonian perturbations for  2D linear   Schr\"odinger equation,
Ann. Math. {\bf 148} (1998) 363--439

\item{[Bo4]}  Bourgain, J., Green function estimates for lattice
Schr\"odinger operators  and applications, preprint

\item{[Bo5]}  Bourgain, J., On Melnikov's persistence problem,
Math. Research Letters {\bf 4} (1997) 445-458

\item{[Br]} Br\'ezis, H., Periodic solutions of nonlinear vibrating strings and duality
principles, Bull. AMS  {\bf 8}  409--426(1983)

\item{[C-W]} Craig,
W. and  Wayne, C.,  Newton's method and periodic
   solutions of nonlinear wave equation,  Commun. Pure. Appl.
   Math.,{\bf  46}   1409--1501 (1993)

\item{[C-Y]} Chierchia, L.  and You, J.: KAM tori for $1$D nonlinear wave equations with periodic boundary conditions.  Comm. Math. Phys. {\bf 211} 497--525 (2000)

\item{[E]} Eliasson, L. H. Perturbations of stable invariant tori
for Hamiltonian systems, Ann. Scula Norm. Sup. Pisa CL Sci. {\bf
15} 115-147 (1988)

 \item{[F-S]} Frohlich,J., Spencer,T. Absence of diffusion in the
 Anderson tight binding model for large disorder or lower energy.
  Commun. Math. Phys.{\bf 88} 151-184 (1983)
%\item{[G-Y]} Geng, J. and You, J., a talk in the conference on
%dynamical systems in Huang Shan, 2005
\item{[G-L-R]} Gohberg, I., Lancaster, P., and Rodman, L.:
   Matrix Polynomials. {\it( Computer Science and Applied Math.,)}
    Academic Press, New York  London (1982)
\item{[Ka]}  Kato, T.:   Perturbation Theory for Linear
    Operators. {\it( Corrected printing of the second edition)}
     Springer-Verlag , Berlin Heidelberg New York (1980)
\item{[K1]} Kuksin, S. B.:  Nearly integrable infinite-dimensional Hamiltonian systems.( Lecture Notes in Math. 1556).
    Springer-Verlag, New York (1993)

\item{[K2]} Kuksin, S. B.:   Elements of a qualitative theory of
Hamiltonian   PDEs, Doc. Math. J. DMV (Extra Volume ICM) {\bf II}
819--829 (1998)

%\item{[K-P]}   Kuksin, S. B. and   P\"oschel, J.: Invariant Cantor manifolds of quasiperiodic    oscillations for a nonlinear Schr\"odinger equation, Ann. Math. {\bf 143}
%      149--179 (1996)
 \item{[L]} Lancaster,P.: Theory of Matrices. Academic Press, New York and London (1969)

\item{[L-S]} Lidskij, B. V. and Shulman, E., Periodic solutions of
the equation $u_{tt}-u_{xx}+u^3=0$, Funct. Anal. Appl. {\bf 22}
332-333 (1988)

\item{[Ph]}  Ph\'ong, V. Q.: The operator equation $AX-XB=C$ with
unbounded operator $A$ and $B$ and related abstract Cauchy
problems. Math. Z. {\bf 208} 567--588 (1991)

\item{[P1]} P\"oschel, J.,  A KAM-theorem for some nonlinear PDEs,
 Ann. Scuola Norm. Sup. Pisa, Cl. Sci., $\quad$IV Ser. 15,
    {\bf 23} 119--148 (1996)

\item{[P2]} P\"oschel, J., Quasi-periodic solutions for nonlinear
wave equation, Commun. Math. Helvetici
    {\bf 71} 269--296 (1996)
\item{[P3]} P\"oschel, J., On the construction of almost periodic
solutions for a nonlinear Schr\"odinger equation, Ergod. Th. \&
Dynam. Syst. {\bf 22} 1537--1549 (2002)

\item{[W]} Wayne, C. E., Periodic and quasi-periodic solutions of nonlinear
    wave equations via KAM theory. Commun. Math. Phys.{\bf 127} 479--528 (1990)

   \item{[Y1]} Yuan, X., Invariant Manifold of Hyperbolic-Elliptic Type for
   Nonlinear Wave Equation, Int. J. Math. Math. Science, {\bf
   18} 1111-1136 (2003)
 \item{[Y2]} Yuan, X., Invariant tori of nonlinear wave
    equations with a given potential.  Preprint. FDIM 2004-11,
    Fudan University (2004), {\it (to appear in DCDS)}.
   \item{[Y3]} Yuan, X., Quasi-periodic solutions of completely
   resonant nonlinear wave equations. Preprint. FDIM 2004-12,
    Fudan University (2004), J. Diff. Equations, (In press).
    %\item{[Y4]} Yuan, X.,  Quasi-periodic solutions of nonlinear
    %Schr\"odinger equations of higher dimension, J. Differential
    %Equations, {\bf 195} 230-242 (2003)
   %item{[Y3]} Yuan, X.: Construction of quasi-periodic breathers
   %via KAM technique.  Commun. Math. Phys.{\bf 226} 61--100 (2002)

\enddocument